\documentclass[11pt]{article}
\usepackage{a4wide}

\usepackage[T1]{fontenc}
\usepackage[utf8]{inputenc}
\usepackage[english]{babel}
\usepackage{amsmath}
\usepackage{amssymb}

\usepackage{color}
\definecolor{grey}{rgb}{.3,.3,.3}

\usepackage{overpic}
\usepackage{babelbib} % Literaturverzeichnis

\usepackage{setspace} % Zeilenabstandspaket 
\usepackage[colorlinks=true, citecolor=black, linkcolor=black, urlcolor=black]{hyperref} %Hyperreflinks auf schwarz stellen

\newcommand\N{{\mathbb N}}
\newcommand\R{{\mathbb R}}
\newcommand\Z{{\mathbb Z}}
\newcommand\Li{L^{\infty}}
\newcommand\T{{\mathcal T}}
\newcommand\K{\mathcal{K}}
\newcommand\Ki{\mathcal{K}_{\infty}}
\newcommand\KL{\mathcal{KL}}
\newcommand\mG{(\Ki \cup \{0\})}
\newcommand\MAF{\mathrm{MAF}}
\newcommand\MAFN{\mathrm{MAF}_{N}}
\newcommand\MAFNN{\mathrm{MAF}_{N}^N}
\newcommand\Gammu{\Gamma_{\mu}}

\newcommand\diam{\operatorname{diam}}
\newcommand\mesh{\operatorname{mesh}}
\newcommand\id{{\operatorname{id}}}

\newcommand\co{{\operatorname{co}}}
\newcommand\1{\operatorname{\textbf{I}}}
\newcommand\2{\operatorname{\textbf{II}}}
\newcommand\3{\operatorname{\textbf{III}}}
\newcommand\4{\operatorname{\textbf{IV}}}
\newcommand\5{\operatorname{\textbf{V}}}

\newtheorem{theorem}{Theorem}[section]
\newtheorem{lemma}[theorem]{Lemma}
\newtheorem{proposition}[theorem]{Proposition}

\newtheorem{definition}[theorem]{Definition}
\newtheorem{remark}[theorem]{Remark}
\newtheorem{assumption}[theorem]{Assumption}

\newenvironment{proof}
{{\par\parskip1.0ex %plus0.5ex minus0.5ex 
\textbf{Proof.}} \setlength{\leftskip}{0.7em}}% Setlength rueckt die Beweise um 1.1 em ein
{$\hfill \square$\parskip1.0ex %plus0.5ex minus0.5ex
\par \vspace{0.5em}}

\begin{document}

%%%%%%%%%%%%%%%%%%%%%%%%%%%%%%%%%%%%%%%%%%%%%%%%%%%%%%%%%%%%%%%%%%%%
%%%	TITLE

\title{Numerical Construction of LISS Lyapunov Functions under a Small Gain Condition}

\author{Roman Geiselhart \and Fabian Wirth\thanks{Email: \href{mailto:roman.geiselhart@mathematik.uni-wuerzburg.de}{\{roman.geiselhart,wirth\}@mathematik.uni-wuerzburg.de}, 
both authors are with \newline Institute for Mathematics, University of W{\"u}rzburg, Emil-Fischer Str. 40, 97074 W{\"u}rzburg, Germany, \newline Tel: +49-931-318 34 36. Fax: +49-931-318 4611}
}

\date{\today}

\maketitle

\begin{abstract}
In the stability analysis of large-scale interconnected
    systems it is frequently desirable to be able to determine a decay
    point of the gain operator, i.e., a point whose image under the
    monotone operator is strictly smaller than the point itself. 
The set of such decay points plays a crucial role in checking, in a semi-global fashion, the local input-to-state stability of an interconnected system and in the numerical construction of a LISS Lyapunov function. 
We provide a homotopy algorithm that computes a
    decay point of a monotone operator. For this purpose we use a fixed
    point algorithm and provide a function whose fixed points correspond
    to decay points of the monotone operator. The advantage to an earlier
    algorithm is demonstrated. Furthermore an example is
    given which shows how to analyze a given perturbed interconnected
    system.
\end{abstract}

\begingroup
\small \textbf{Keywords:} homotopy algorithm, monotone operator, LISS
Lyapunov function, interconnected system, small gain condition
\endgroup

\medskip

In recent years large-scale systems have received renewed attention with
applications in formation control, logistics, consensus dynamics,
networked control systems and further applications. While stability
conditions for such large-scale systems have already been
studied in \cite{moylan1978stability,Sil79,Vid81} based on linear gains
and Lyapunov techniques, nonlinear approaches are more recent. The
groundbreaking concept that has proven fruitful is the notion of
input-to-state stability (ISS) as introduced in \cite{Son89}.

For large-scale nonlinear systems it may be difficult to prove ISS
directly, but
if a large-scale system is defined through the interconnection of a number
of smaller components, which are ISS, then there exist small gain type
conditions guaranteeing the ISS property for the interconnected
system. For the case of two subsystems this result was obtained in
\cite{JPT94,JMW96} both in a trajectory based as well as a Lyapunov
formulation. Recently, there has been a substantial effort to extend these
results to the case of a greater number of subsystems, see
\cite{DRW07,DRW09,ito2009necessary,liu2009lyapunov,rueffer2010small,DashIto2011}.
It is the purpose of this paper to provide numerical methods that make some of the available results applicable for practical problems.

The general setting is here to consider a number of systems that are
input-to-state stable with respect to external and internal inputs.  The
effect of the subsystems, described by comparison functions, is collected
in the \emph{gain matrix} $\Gamma$. The special structure of the
interconnected system now leads to a monotone operator $\Gammu$ on the
positive orthant $\R^N_+$. 
So-called monotone aggregation functions can be used to formulate the effect of several inputs on a system in a general manner.
Standard examples of such functions are summation and
maximization, but in \cite{DRW09} some examples are provided that also
other types of aggregation functions may be useful depending on the system
under consideration. We would like to point out that the particular
relation of the maximization and summation formulation of small gain
conditions is analyzed in \cite{DashKosm11}.
In \cite{Karafyllis:2010fk} the authors study interconnections where small gain conditions are satisfied after certain transient periods and derive stability results.

Many available small gain results state that input-to-state stability for
the overall system follows from the existence of a so-called $\Omega$-path
with respect to $\Gammu$,
\cite{DashIto2011,DashKosm11,DR09,DRW07,DRW09,Ruffer:2009:Monotone-inequalities-dynamical-systems-::}.
Furthermore an ISS Lyapunov function for the interconnected system can be
constructed using this path and the ISS Lyapunov functions of the
subsystems.  Note that also for other small gain type formulations as the
cycle condition in the maximization case or the spectral radius condition
in the linear summation case, it may be seen that these conditions can be
equivalently formulated in terms of $\Omega$-paths.

In \cite{DRW09} the construction of an $\Omega$-path is described. The
crucial ingredient that usually cannot be obtained in a straightforward
manner is a \emph{decay point} of $\Gammu$, that is a point $s \in \R^N_+$
for which $\Gammu(s) \ll s$ in the order induced by the cone
$\R^N_+$. Once such a point is found there are straightforward numerical
procedures for the construction of Lyapunov functions or for checking the
ISS property. There are two particular cases in which a straightforward
way is known to compute decay points: $(i)$ If the gains are linear and
summation is used, then the problem becomes one of checking whether the
spectral radius of $\Gamma$ is below 1 and finding an appropriate
eigenvector. As $\Gamma$ is nonnegative this problem is particularly easy
and well studied; $(ii)$ If the maximization formulation of ISS is used
then a very nice observation of \cite{KJ09} is that, provided a
small gain-condition holds, for $s^*:=
\max\{s,\Gamma(s),\ldots,\Gammu^{N-1}(s)\}$, we have $\Gamma(s^*)\leq s^*$
for all $s \in \R^N_+$, which is almost a decay point. The methods
presented in this paper are suitable for the cases that the problem at
hand is not within one of the two classes described above.
In this paper we provide numerical procedures for computing such points and thus also for local $\Omega$-paths. 
We call the approach \emph{semi-global} because it does not require a priori restrictions. In particular, if a small gain condition is  satisfied globally, then the design variables of the algorithm can in principle be chosen so that the numerically guaranteed region of stability is arbitrarily large.

As we compute a decay point numerically the overall construction of
Lyapunov functions as well as the verification of the ISS property is only
performed locally. Indeed, the approach relies on local results of the
small gain type. Local small gain results have been considered in
\cite{BourColl95,Bour2000,JianLin2004} in an input-output operator
context, resp. for discrete-time systems. In \cite{DR09} local ISS (LISS)
definitions and local small gain theorems within the framework considered
here. In this work the knowledge of a decay point leads to the local
input-to-state stability of the interconnected system and to the
construction of a LISS Lyapunov function.

The algorithm developed here, that computes a decay point for a given monotone operator
$\Gammu$, is a particular \emph{simplicial fixed point algorithm} (SFP-algorithm) customized in such a way that we obtain a decay point of
$\Gammu$. To ensure the convergence of the SFP-algorithm we require
irreducibility of the gain matrix $\Gamma$. This is no significant restriction because by
standard graph theoretic algorithms the irreducible components of the
system can be obtained efficiently, \cite{Tar72}. 

The paper is organized as follows. In Section~\ref{sec:Pre} we provide the
necessary notions and a short introduction to comparison functions and
graphs. In Section~\ref{sec:ISS} we recall the Lyapunov formulation of ISS
for interconnected systems, give a local small gain theorem and outline
the construction of a LISS Lyapunov function for the overall
systems. Section~\ref{sec:Merrill} contains the main results of this
paper. First we recall some facts about homotopy algorithms and introduce
the SFP-algorithm where we mainly follow the book of \cite{Yang}. In
subsection \ref{subsec:Anpassung} we state some sufficient conditions on
$\Gammu$ and prove that the SFP-algorithm converges to a decay point of
$\Gammu$. At the end of this section some improvements of the algorithm
are discussed. We conclude this work in Section~\ref{sec:Exa} where we
discuss two examples. The first one shows that this new algorithm improves
on an earlier algorithm that is due to a homotopy algorithm of Eaves \cite{Eaves}
(cf. \cite{RufferWirth:2010:Stability-verification-for-monotone-syst:})
where we revisit a nonlinear example from
\cite{RufferDowerIto:2010:Applicable-comparison-principles-in-larg:}. In
the second example we use our techniques to show numerically that a
particular perturbed interconnected system is LISS.

%%%%%%%%%%%%%%%%%%%%%%%%%%%%%%%%%%%%%%%%%%%%%%%%%%%%%%%%%%%%%%%%%%%%
%%% 	PRELIMINARIES
\section{Preliminaries}\label{sec:Pre}

\subsection{Notation and conventions}
\label{subsec:Grundlagen-Notation}

Let $\R$ denote the field of real numbers, $\R_+$ the set of nonnegative
real numbers, and $\R^N$ (resp. $\R^N_+$) the vector space of
(nonnegative) real column vectors of length $N$. Then $\R^N$ induces a
partial order for vectors $v,w \in \R^N$. We denote $v \geq w \iff
v_{i}\geq w_{i}$, $v > w \iff v\geq w$ and $ v \neq w$, $v \gg w \iff
v_{i} > w_{i}$, each for $i=1, \ldots, N$, where $v_i$ denotes the
$i^{th}$ component of the vector $v$.  Let $v,w \in \R^N_{+}$ be
given. Then we define the \emph{order intervals} $[v,w] := \{ x \in
\R^N_{+} : v \leq x \leq w \}$ if $v \leq w$, $(v,w) := \{ x \in \R^N_{+}
: v \ll x \ll w \}$ if $v \ll w$, and analogously the order intervals
$(v,w]$ and $[v,w)$.  For $x \in \R^N$ we use the Euclidean norm $\|x\| =
\sqrt{\sum_{i=1}^N |x_{i}|^2}$.  The space of measurable and essentially
bounded functions is denoted by $\Li=\Li([0,\infty); \R^M)$ with norm $\|
\cdot\|_\infty$.

\subsection{Comparison functions and induced monotone operators}
\label{subsec:Grundlagen-Vergleichsfunktionen}

To state the stability definitions that we are interested in, three sets
of comparison functions are used.  We call a function $\alpha: \R_{+}
\rightarrow \R_{+}$ a \emph{function of class $\K$}, if it is strictly
increasing, continuous, and satisfies $\alpha(0)=0$. If $\alpha \in \K$ is
unbounded, it is said to be of class $\Ki$.  A function $\beta:\R_{+}
\times \R_{+} \rightarrow \R_{+}$ is called a \emph{function of class
  $\KL$}, if it is of class $\Ki$ in the first argument and strictly
decreasing to zero in the second argument.  It is easy to see that if
$\rho \in \Ki$, then its inverse $\rho^{-1}: \R_{+} \rightarrow \R_{+}$
exists and is also of class $\Ki$.

To formulate general small gain conditions we need the following
definition, see \cite{DRW09}.
\begin{definition}\label{Def:MAF} 
    A continuous function $\mu: \R_{+}^N \rightarrow \R_{+}$ is called a
    \emph{monotone aggregation function}, if the following properties
    hold:
\begin{enumerate} \item positivity: $\mu(s) \geq 0$ for all $ s \in \R^N_{+}$ and $\mu(s)=0$, if and only if $s=0$;
\item strict increase: $ \mu(s) < \mu (t)$, if $s\ll t$;
\item unboundedness: $\mu(s) \rightarrow \infty$, if $\|s \| \rightarrow \infty$.
\end{enumerate}
The space of monotone aggregation functions is denoted by $\MAFN$. 
\end{definition}

The properties in Definition \ref{Def:MAF} can be extended to vectors in the sense that $\mu = (\mu_{1}, \ldots, \mu_{N})^\top \in \MAFN^N$, $\mu_{i}\in \MAFN, i=1, \ldots, N$, defines a mapping from $\R^{N \times N}$ to $\R^N$ by
$ \mu(A)_{i}= \mu_{i}(a_{i1}, \ldots, a_{iN})$ for $A= (a_{ij})_{i,j=1}^N \in \R^{N \times N}_{+}$.\\
We want to generalize this to matrices of the form $\Gamma=(\gamma_{ij})_{i,j=1}^N \in (\Ki \cup \{0\})^{N \times N}$, where $0$ denotes the zero function. This leads to an operator $\Gammu: \R^N_{+} \rightarrow \R^N_{+}$ defined by
\begin{equation}\label{eq:GammuOp} \Gammu(s) :=(\mu \circ \Gamma)(s):= \left(\begin{array}{c} \mu_{1}(\gamma_{11}(s_1), \ldots, \gamma_{1N}(s_{N}) )\\\vdots \\ \mu_{N}(\gamma_{N1}(s_1), \ldots, \gamma_{NN}(s_{N}) )\end{array}\right) \in \R^N_{+} \quad \text{ for } s \in \R^N_{+}. \end{equation}
For the $k$ times composition of this operator we write $\Gammu^k$. We call $\Gammu$
\begin{enumerate}
\item \emph{monotone}, if  $\Gammu(v) \leq \Gammu(w)$ for all $v,w \in \R^N_{+}$ with $v \leq w$;
\item \emph{strictly increasing}, if $\Gammu(v) \ll \Gammu(w)$ for all $v,w \in \R^N_{+}$ with $v \ll w$.
\end{enumerate}

\begin{remark}\label{remark:ground-Gammu}
Note that if $\Gamma \in (\Ki \cup \{0\})^{N\times N}$ and $\mu \in \MAFNN$, then $\Gammu$ is monotone and satisfies $\Gammu(0)=0$.
\end{remark}

The next definition is fundamental in the following.

\begin{definition}\label{Def:AbMeng} For a given function $T: \R_{+}^N \rightarrow \R_{+}^N$ we define the \emph{set of decay} $\Omega$ by $$ \Omega(T):= \left\{ s \in \R_{+}^N: T(s)\ll s \right\}.$$
For short we just write $\Omega$, if the reference to $T$ is clear from the context. Points in  $\Omega$ are called \emph{decay points}. \end{definition}

\subsection{Graphs and matrices}\label{subsec:Grundlagen-Matrizen}

A directed graph $G(V,E)$ consists of a finite set of vertices $V$ and a set of edges $E \subset V \times V$. If $G(V,E)$ consists of $N$ vertices, then we may identify $V=\{1, \ldots, N\}$. So if $(j,i) \in E$, then there is an edge from $j$ to $i$. The \emph{adjacency matrix} $A_{G}=(a_{ij})$ of this graph is defined by $a_{ij} = 1$, if $(j,i) \in E$ and $a_{ij}=0$ else. 
We call the graph $G(V,E)$ \emph{strongly connected}, if for each pair $(i,j)$ there exists a \emph{path} 
$(e_{i_{0},i_{1}}, e_{i_{1},i_{2}}, \ldots, e_{i_{k-1},i_{k}})$ with $i=i_{0}, j=i_{k}$ such that $e_{i_{l-1}, i_{l}} \in E$ for all $i=1, \ldots, k$. 
It is well known that the graph $G(V,E)$ is strongly connected, if and only if the adjacency matrix $A_{G}$ is irreducible, i.e., there exists no permutation matrix $P$ such that $$ A= P^T  \left(\begin{array}{cc}B & C \\0 & D\end{array}\right) P$$ for suitable, square matrices $B$ and $D$. These definitions can be carried over to matrices $\Gamma \in \mG^{N \times N}$. To this end we define the matrix $A_{\Gamma}=(a_{ij})_{i,j=1}^N$ by $a_{ij} = 1$, if $\gamma_{ij} \in \Ki$, and $a_{ij}=0$, if $\gamma_{ij}\equiv 0$.  We call $\Gamma$ irreducible, if the matrix $A_{\Gamma}$ is.

%%%%%%%%%%%%%%%%%%%%%%%%%%%%%%%%%%%%%%%%%%%%%%%%%%%%%%%%%%%%%%%%%%%%
%%%	INPUT-TO-STATE STABILITY

\section{Input-to-state stability and small gain theorems}\label{sec:ISS}

Consider the control system
\begin{equation}\label{eq:ISS-nlgdglu} \dot x(t) = f(x(t), u(t)), \hspace{2em} t \in \R_{+}, \end{equation}
where $u \in \R^m$ is the \emph{input} and $x \in \R^n$ is the
\emph{state}. We assume that $f: \R^n \times \R^m \rightarrow \R^n$ is
continuous and \emph{locally Lipschitz in $x$ uniformly for $u$ in
  compacts}; by this we mean that for every compact $K_1\subset\R^n$ and
compact subset $K_2\subset\R^m$ there is some constant $c>0$ such that
$\|f(x,u)- f(z,u)\| \leq c \|x-z\|$ for all $x,z \in K_{1}$ and all $u \in
K_{2}$.  Further we assume $f(0,0)=0$ and all solutions can be extended to
$[0, \infty)$.

\begin{definition}\label{def:ISS-Lya} 
    Consider the system $(\ref{eq:ISS-nlgdglu})$ and let $V:\R^n
    \rightarrow \R_{+}$ be continuous and locally Lipschitz continuous
     on $\R^n\backslash \{0\}$. Then $V$ is called an \emph{ISS
      Lyapunov function}, if there exist $\alpha_{1}, \alpha_{2}\in \Ki$
    such that for all $x \in \R^n$,
    \begin{equation}\label{eq:ISS-Lyaprop1} \alpha_{1}(\|x\|) \leq V(x)
        \leq \alpha_{2}(\|x\|), \end{equation}
    and if there exist $\gamma \in \K$ and a positive definite function 
    $\alpha_{3}$ 
    such that for all $u \in \R^m$ and almost all $x \in \R^n$, 
    \begin{equation}\label{eq:ISS-Lyaprop2} V(x)\geq \gamma(\|u\|) \quad
        \Rightarrow \quad \nabla V(x)f(x,u) \leq
        -\alpha_{3}(\|x\|). \end{equation}
\end{definition}
Note, that we only assume Lipschitz continuity of the ISS Lyapunov
function $V$. By Rade\-macher's Theorem, see e.g. \cite{Rademacher}, this
implies that $V$ is differentiable almost everywhere and we consider the
decay condition \eqref{eq:ISS-Lyaprop2} only at points where $V$ is differentiable. An equivalent
formulation can be given in terms of Clarke subdifferentials but we
refrain from doing so, since this will play no further role in the paper,
see also \cite{DRW07,DRW09}.

System $(\ref{eq:ISS-nlgdglu})$ is called input-to-state stable (ISS), if
it has an ISS Lyapunov function. There is another, trajectory-based
definition of ISS which is equivalent to the existence of an ISS Lyapunov
function (cf. \cite{SW95} for smooth Lyapunov functions and \cite[Theorem
2.3]{DRW09} for continuous and locally Lipschitz continuous functions).

Now we want to generalize this stability definition to networks. Let
$N\in\N$ and consider the $N$ interconnected systems given by
\begin{equation}\label{eq:ISS-ics} \begin{array}{ccc}\dot x_1 & = &
        f_1(x_1,\dots,x_N,u) \\ & \vdots & \\\dot x_N & = &
        f_N(x_1,\dots,x_N,u)\end{array}. \end{equation}
Assume that $x_{i} \in \R^{n_{i}}, u \in \R^m$ and  the functions $f_{i}: \R^{\sum_{j=1}^N n_{j}+m} \rightarrow \R^{n_{i}}$ are continuous and locally Lipschitz in $x=(x_{1}^\top, \ldots, x_{N}^\top)^\top$ uniformly for $u$ in compacts. Let $x_{i}$ denote the state of the $i^{th}$ subsystem and assume $u$ as an external control variable. 
Without loss of generality we may assume to have the same input for all systems, since we may consider $u$ as partitioned $u=(u_{1}^\top, \ldots, u_{N}^\top)^\top$, such that each $u_i$ is the input for subsystem $i$ only.  Then each $f_{i}$ is of the form $f_{i}(\ldots, u) = \tilde f_{i}(\ldots, \pi_{i}(u)) = \tilde f_{i}(\ldots, u_{i})$ with a projection $\pi_{i}$.

If we consider individual systems, we treat the state $x_j, j \neq i$, as an independent input for $x_i$. Assume that for each subsystem $i\in \{1, \ldots, N\}$ there exists a continuous and locally Lipschitz continuous function $V_{i}:\R^{n_{i}} \rightarrow \R_{+}$ such that for suitable $\alpha_{1i}, \alpha_{2i}\in \Ki$
\begin{equation}\label{eq:ISS-icsLyaprop1} \alpha_{1i}(\|x_{i}\|) \leq V_{i}(x_{i}) \leq \alpha_{2i}(\|x_{i}\|) \quad \text{ for all } x_{i} \in \R^{n_{i}}. \end{equation}
We call the function $V_{i}$ an \emph{ISS Lyapunov function} for the subsystem $i$, if there exist $\mu_{i} \in \MAF_{N+1}$, $\gamma_{ij} \in \Ki\cup \{0\}, 
j \neq i, \gamma_{iu}\in \K \cup \{0\}$ and a positive definite function $\alpha_{i}$ such that
\begin{eqnarray}\label{eq:ISS-icsLyaprop2} V_{i}(x_{i}) &\geq& \mu_{i} \big( \gamma_{i1}(V_{1}(x_{1})), \ldots,  \gamma_{iN}(V_{N}(x_{N})), \gamma_{iu}(\|u\|)  \big) \nonumber \\
&& \Rightarrow \ \nabla V_{i}(x_{i})f_{i}(x,u) \leq - \alpha_{i}(\|x_{i}\|). \end{eqnarray}
%\setlength{\arraycolsep}{5pt}
%\begin{equation}\label{eq:ISS-icsLyaprop2} V_{i}(x_{i}) \geq \mu_{i} \big( \gamma_{i1}(V_{1}(x_{1})), \ldots,  \gamma_{iN}(V_{N}(x_{N})), \gamma_{iu}(\|u\|)  \big)  \ \Rightarrow \ \nabla V_{i}(x_{i})f_{i}(x,u) \leq - \alpha_{i}(\|x_{i}\|). \end{equation}
The functions $\gamma_{ij}$ and $\gamma_{iu}$ are called \emph{ISS Lyapunov gains}. We distinguish between the \emph{internal inputs} $x_{j}$ and the \emph{external input} $u$ of the $i^{th}$ subsystem. These gains indicate the influence of the inputs on the state. This is why we set $\gamma_{ij}\equiv 0$, if $f_{i}$ does not depend on $x_{j}$ and we collect the internal inputs into the \emph{gain matrix} $\Gamma:=(\gamma_{ij})_{i,j=1}^N$. Note that $\Gamma$ and the $\mu_i$ define a monotone operator $\Gammu: \R_{+}^{N} \rightarrow \R_{+}^{N}$ as in (\ref{eq:GammuOp})  
(cf. Remark \ref{remark:ground-Gammu}).

\subsection{A local small gain theorem}\label{subsec:ISS-sgt}

In this section we assume that the interconnected system (\ref{eq:ISS-ics}) satisfies an ISS condition of the form (\ref{eq:ISS-icsLyaprop2}) for ISS Lyapunov functions $V_i$, $i=1, \ldots, N$. Denote the corresponding gain operator by $\Gammu$ as in (\ref{eq:GammuOp}). We assume that $\Gamma$ is irreducible, so that $\Gammu$ is strictly increasing (cf. \cite[Lemma 2.7]{Ruffer:2009:Monotone-inequalities-dynamical-systems-::}). A local ISS Lyapunov function for the overall system given by
\begin{equation}\label{eq:ISS-ics2} \dot x = f(x,u) \end{equation} 
and $x=(x_{1}^\top, \ldots, x_{N}^\top)^\top, f = (f_{1}^\top, \ldots, f_{N}^\top)^\top$ may now be constructed as follows.\\
Assume there exists a $w\gg 0$ with
\begin{equation} \label{eq:decayw}
 \Gammu(w) \ll w.
\end{equation}
Then the sequence $\Gammu^k(w)$, $k=1,2, \ldots$ is strictly decreasing and so $\lim_{k\rightarrow \infty} \Gammu^k(w)$ exists. If
\begin{equation} \label{eq:lim=0}
 \lim_{k \rightarrow \infty} \Gammu^k(w) =0, 
\end{equation}
then we define the linear interpolation of the points $\{ \Gammu^k(W) \}_{k\in \N}$ by $\sigma:[0,1]\rightarrow \R^N_+$:
\begin{equation}\label{eq:lininterpol}
\sigma(r) = \left\{\begin{array}{ll}
 0, & \text{if } r=0 \\ 
 (k^2+k)\left( [\tfrac1k -r]\Gammu^k(w) + [r-\tfrac1{k+1} ]\Gammu^{k-1}(w) \right), & \text{if } r \in (\tfrac{1}{k+1}, \tfrac{1}{k}], \ k \in \N.
\end{array} \right.
\end{equation}

Note that $\sigma$ is continuous on $[0,1]$ by (\ref{eq:lim=0}) and strictly increasing in all component functions as $\Gammu$ is assumed to be irreducible. With this construction local ISS Lyapunov functions can be constructed using the following summary of existing results (cf. \cite[Theorem 5.5]{DR09}).

\begin{theorem}\label{theo:LISS}
 Assume that system $(\ref{eq:ISS-ics})$ satisfies ISS conditions of the form $(\ref{eq:ISS-icsLyaprop2})$ for all $i=1,\ldots,N$, and that the gain matrix $\Gamma$ is irreducible. If there exists an $w\gg0$ so that $(\ref{eq:decayw})$ and $(\ref{eq:lim=0})$ hold, then a local ISS Lyapunov function for the overall system $(\ref{eq:ISS-ics2})$ is given by
\begin{equation}\label{eq:LISS-Vnetwork} V(x)= \max_{i=1, \ldots, N} \sigma_{i}^{-1}(V_{i}(x_{i})).
\end{equation}
In particular, the implication
\begin{equation}
\label{eq:smallgaincondloc}
 V(x) \geq \gamma (\|u\|) \quad \Rightarrow  \quad \nabla V(x) \cdot f(x,u) \leq - \alpha(V(x)) \end{equation}
holds locally with $\gamma \in \Ki$ given by \cite[Proposition 4.3]{DR09}. 
\end{theorem}

\begin{remark}
 \begin{enumerate}
   \item By ``local'' we mean ``in an open neighborhood of the origin
     $(x^*,u^*)=(0,0)$''. In particular, \cite[Theorem 5.5]{DR09} shows
     that the assured domain of stability increases with the choice of
     $w$. In particular, if the small gain condition holds globally the domain where
     \eqref{eq:smallgaincondloc} holds can be made arbitrarily large.
  \item By $(\ref{eq:decayw})$ and $(\ref{eq:lim=0})$ we have the \emph{small gain condition} $\Gammu(s) \not \geq s$ for all $s \in [0,w]$.
  \item Note that $\sigma(r) \in \Omega(\Gammu)$ for all $r \in [0,1]$ and
    $\sigma$ belongs to the class of $\Omega$-paths (cf. \cite[Definition
    5.1]{DRW09}).
 \end{enumerate}
\end{remark}

\begin{remark}
    (i) Theorem \ref{theo:LISS} is the starting point for our numerical
    considerations. If we find the decay point $w$, then the problem of
    constructing Lyapunov functions or checking small gain conditions
    becomes easy. In the remainder of the paper we concentrate on
    giving numerically tractable solutions to this problem.
 
    (ii) In the linear case with $\mu=\Sigma$ we have $\Gamma_\Sigma(s)= \Gamma
    s$ with $\Gamma \in \R_+^{N\times N}$. Here the existence of a decay
    point $w\gg0$ with $\Gamma w \ll w$ is equivalent to the spectral
    radius of $\Gamma$ being less than one, i.e., $1>\rho(\Gamma)=
    \{|\lambda| : \lambda \text{ is an eigenvalue of } \Gamma \}$
    (cf. \cite[Lemma
    1.1]{Ruffer:2009:Monotone-inequalities-dynamical-systems-::}). So
    finding a decay point is just an eigenvalue problem. This is why we
    assume $\Gammu$ to be nonlinear.
\end{remark}

%%%%%%%%%%%%%%%%%%%%%%%%%%%%%%%%%%%%%%%%%%%%%%%%%%%%%%%%%%%%%%%%%%%%
%%% 	MERRILL'S ALGORITHM 

\section{A homotopy algorithm for computing a decay point $\mathbf{w \in \Omega(\Gammu)}$}\label{sec:Merrill}

In this section we want to develop an algorithm that computes a decay
point $w \in \Omega(\Gammu)$ for a given continuous and monotone operator
$\Gammu: \R^N_{+} \rightarrow \R^N_{+}$.  We know that such a point exists
for any norm, if the \emph{small gain condition}
\begin{equation}\label{eq:Merrill:SGB} \Gammu(s) \not \geq s \quad \text{for all } s \in \R^N_{+}\backslash \{0\}\end{equation} 
is satisfied (cf. \cite[Proposition 5.3]{DRW07}). 

To find such a point we will extend a homotopy algorithm that was also
used by Merrill \cite{Mer} to compute fixed points of upper-semicontinuous
(u.s.c.) point-to-set mappings. Note that since a continuous single-valued
function is in particular an u.s.c. point-to-set mapping our problem falls
in the class of problems that can be treated by homotopy
algorithms. However, Merrill's condition introduced in \cite{Mer} is not
sufficient for convergence in our case, as the domain of the mapping is
only the nonnegative orthant.  The idea to the design of a 
convergent algorithm is to construct a function $\phi: \R^N_+ \rightarrow
\R^N$, which has the property that its fixed points are decay points of
$\Gammu$, and to show that the homotopy algorithm will converge to
approximate fixed points of $\phi$, which are also decay points 
of $\Gammu$. This algorithm is semi-global since by choosing design variables appropriately we end up in a decay point with arbitrarily large norm.

In Section~\ref{subsec:Grundlagen-Simplizes} we present the triangulation
we need for the computation of fixed points. Before introducing the
homotopy algorithm in Section~\ref{subsec:Algo} we first provide some
facts about homotopy algorithms in Section~\ref{subsec:Homotopie}. In
Section~\ref{subsec:Algo} we mainly follow the book of Yang
\cite[Section 4.3]{Yang}. In Section~\ref{subsec:Anpassung} we will
give the function $\phi$ mentioned above and prove the convergence of the
SFP-algorithm. Finally in Section~\ref{subsec:ImproveMerrill} we give
approaches for improving the algorithm and further give suggestions for
the choice of the design variables used in the mapping $\phi$.

\subsection{Simplices and triangulations}\label{subsec:Grundlagen-Simplizes}

We briefly recall facts about covering convex sets by triangulations. A
set $C \subset \R^N$ is called \emph{convex}, if for all $a,b\in C$ it
holds $S_{a,b}:= \{\lambda a + (1- \lambda)b : \lambda \in [0,1]\} \subset
C$. The {\em convex hull} 
$\co\{ M \}$ of a set $M\subset \R^N$ is the
smallest convex set containing $M$. If $M$ is finite, we also say that
$\co\{ M \}$ is {\em spanned} by $M$ and denote this by $\langle \{v^i \in M \} \rangle$. The \emph{dimension} of a convex set is equal
to the dimension of the smallest affine subspace $U \subset \R^N$
containing $C$.

\begin{definition}\label{Def:Simplex} An $N$-\emph{simplex} $S$ is an
    $N$-dimensional, convex polytope spanned by $N+1$
    vectors $v^{1}, \ldots, v^{N+1}$ in $\R^{M}$, $M \geq N$, i.e.,
$$ S=\langle  v^1, \ldots, v^{N+1} \rangle.$$ A
    \emph{subsimplex} $\varsigma$ of $S$ is a simplex spanned by a subset of
    the set of vertices of $S$, i.e., $\varsigma=\{v^{i} : i \in I_{\varsigma}\}$
    with $I_{\varsigma} \subset \{1, \ldots, N+1\}$. Zero-dimensional
    subsimplices are just the vertices of the simplex, one-dimensional
    subsimplices are called \emph{edges} between the vertices and
    $(N-1)$-subsimplices are called \emph{facets}.
    The subsimplex $S(j)=\langle \{v^i : i \neq j\} \rangle$ is called the
    \emph{facet opposite} $v^j$.
\end{definition}

Clearly, since any $N$-simplex is $N$-dimensional, $N$ of the $N+1$
vertices are linearly independent and it holds $v^i\neq v^j$ for $i\neq
j$. Simplices can be used to cover convex sets in $\R^N$ as follows.

\begin{definition}\label{def:triang} Let $C$ be an $m$-dimensional convex
    set in $\R^N$. A set $\T$ of $m$-simplices is called a
    \em{triangulation} of $C$, if
\begin{enumerate} \item $C$ is the union of all simplices in $\T$;
\item for any $\eta_{1}, \eta_{2} \in \T$, $\eta_1\neq \eta_2$, the intersection $\eta_{1} \cap \eta_{2}$ is either the empty set or a common facet of both;
\item every $x\in C$ has an open neighborhood intersecting only a finite number of $\eta\in\T$. 
\end{enumerate}\end{definition}

By $\T^k$ we denote the set of all $k$-subsimplices of $\T$. It is easy to see that $\T^N = \T$ and $\T^0$ describes the set of the vertices of the simplices in $\T$. To distinguish simplices, or triangulations, we introduce the \emph{diameter of a simplex} $\eta \in \T$ by
$$ \diam(\eta) = \max \{ \| x-y\| :  x,y \in \eta \}$$
and the \emph{mesh size of a triangulation} $\T$ by
$$\mesh(\T) = \sup \{ \diam(\eta) :  \eta \in \T \}. $$

\begin{figure}[ht]
\centering
\begin{overpic}[scale=0.3505]{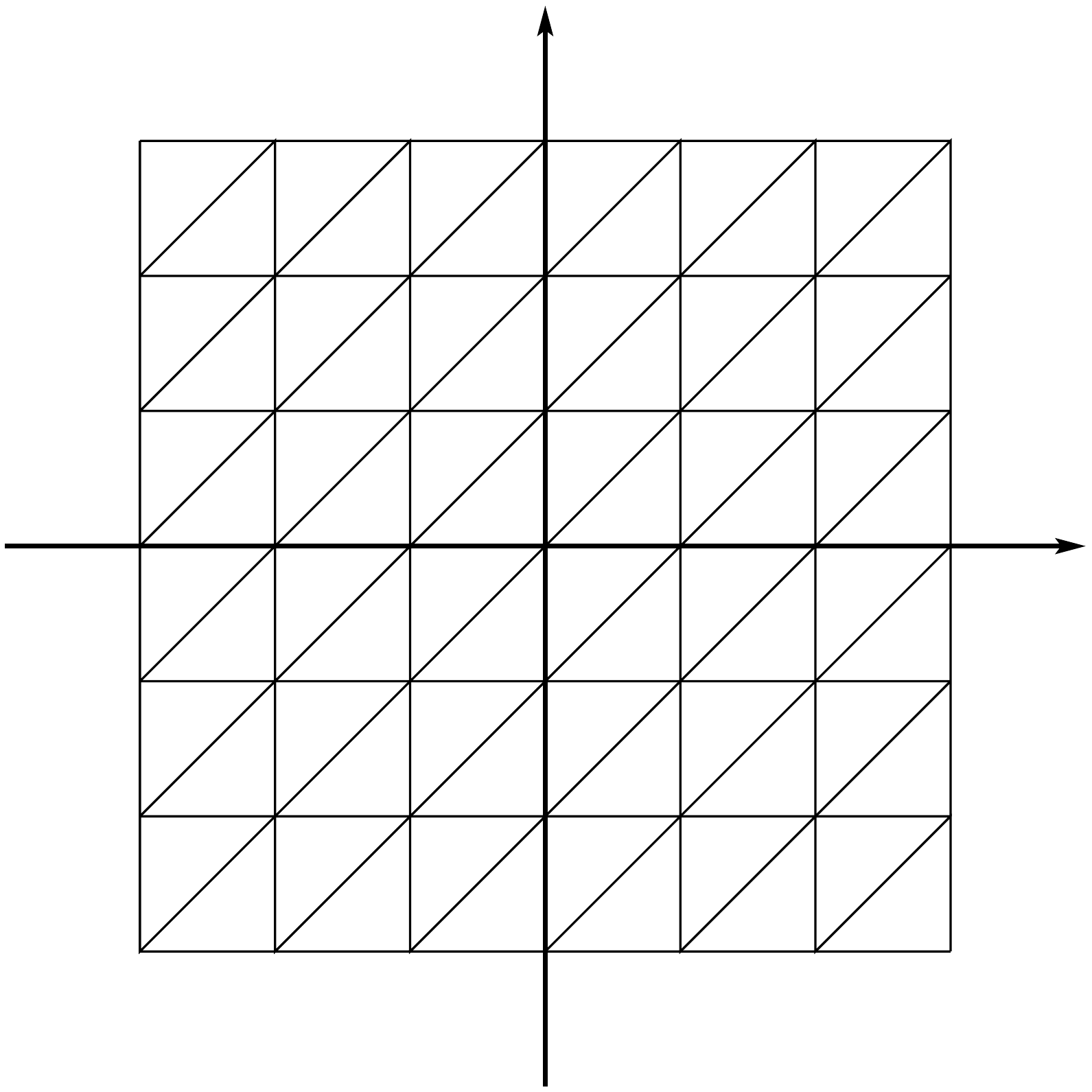}
\put(100.3,49.3){${\scriptstyle x_{1}}$} 
\put(48,100.5){${\scriptstyle x_{2}}$}
\end{overpic}
\caption{Illustration of the $K_{1}$-triangulation for $N=2$}
\label{fig:KEinsTriang}
\end{figure}

There is one special triangulation of $\R^N$, which will be used to
compute decay points. Let $e_{i}$ denote the $i^{th}$ unit vector in
$\R^{N}$. The $K_1$-triangulation is defined as the set of all
$N$-simplices with vertices $x^1, \ldots, x^{N+1}$ such that
$$x^1 \in \Z^N \text{ and } x^{i+1} = x^i+ e_{\pi_{N}(i)}\text{ for all } i \in \{1, \ldots, N\},$$
where $\pi_{N} = (\pi_{N}(1), \ldots, \pi_{N}(N))$ is a permutation of the elements of the set $\{1, \ldots, N\}$. We denote these simplices by $\eta(x^1, \pi_{N})$. %
See \cite[Theorem 1.4.8]{Yang} for a proof that $K_1$ is a triangulation
in the sense of Definition~\ref{def:triang}. An illustration of this
triangulation is given in Figure \ref{fig:KEinsTriang}.

Defining $\delta C= \{ \delta x : x \in C\}$ for $C \subset \R^N$, $\delta>0$, and $\delta F = \{ \delta C : C \in F \}$ for a family $F$ of subsets of $\R^N$ we obtain that if $\T$ is a triangulation of $C$ and $\delta >0$, then $\delta \T$ is a triangulation of $\delta C$. In this way we get the $\delta K_{1}$-triangulation of $\R^N$ for which $\mesh(\delta K_{1})=\delta \sqrt{N}$ for $\delta >0$.\\
Let $\T$ be a triangulation of $\R^N \times [0,1]$ with the restriction $\T^{0} \subseteq \R^N \times \{0,1\}$, i.e., the vertices only lie in $\R^N \times \{0,1\}$. Then we call this triangulation \emph{two-layered}.
Let $\tilde K_{1}$ denote the restriction of the  $K_{1}$-triangulation of $\R^{N+1}$ to $\R^N \times [0,1]$. Then $\tilde K_{1}$ is two-layered. Further define the $(N+1)\times(N+1)$-matrix $P=[\delta e_{1}, \ldots, \delta e_{N}, e_{N+1}]$ for given  $\delta>0$. Define
$$\tilde K_{1}(\delta) = \{ \langle Py^1, \ldots, Py^{N+2} \rangle : \langle y^1, \ldots, y^{N+2}\rangle \in \tilde K_{1} \},$$
then $\tilde K_{1}(\delta)$ is a two-layered triangulation of $\R^N \times [0,1]$.

\subsection{Some facts about homotopy algorithms}\label{subsec:Homotopie}

In this section we want to provide the basic principles of homotopy algorithms. 

\begin{definition} Let $f,g : C\rightarrow D$ be two continuous mappings from the topological space $C$ to the topological space $D$. We call $f,g$ \emph{homotopic}, if there exists a continuous mapping $\vartheta : C\times [0,1] \rightarrow D$, $(s,t)\mapsto \vartheta(s,t)$ with $\vartheta(s,0) =f(s)$ and $\vartheta(s,1) = g(s)$ for all $s \in C.$ 
We call $\vartheta$ the \emph{homotopy} from $f$ to $g$.\end{definition} 

Let $C$ be a nonempty, compact and convex subset of $\R^N$ and assume that $f: C \rightarrow C$ is continuous. Then it follows by Kakutani's fixed point theorem (cf. \cite[p.174]{Ber97}) that there exists at least one fixed point of $f$. To determine any fixed point we use the idea of the classical homotopy. Define the continuous mapping $f_{t}: C \rightarrow C$ by 
$$f_{t}(s):= (1-t)s_{0}+tf(s), \qquad t \in [0,1]$$ 
with $s_{0} \in C$. Then by a further application of Kakutani's fixed point theorem, there exists a fixed point of $f_t$ for every $t \in [0,1]$. We start with the constant mapping $f_0(s)=s_0$ and its fixed point $s_0$. Assume that $t_k\rightarrow 1$ for $k \rightarrow \infty$, then the sequence of functions $(f_{t_k}(\cdot))_{k \in \N}$ converges even uniformly to $f_{1}(\cdot)=f(\cdot)$. 
Now one can show that the cluster points of the set of fixed points $s_{t_{k}}$ of $f_{t_{k}}$ are just the fixed points of $f$. Note that in this approach we have to extend the dimension of this problem, i.e., we now work in the space $C \times [0,1]$.\\
The numerical procedure for nonempty, compact and convex $C \subseteq
\R^N$ is the following. We decompose the space $C \times [0,1]$ in
simplices using a suitable triangulation $\T$. Under certain conditions
there exists a path in this triangulation from an $N$-simplex $\tau^0 \in
C \times \{0\}$ to an $N$-simplex $\tau^* \in C \times \{1\}$ which yields
an approximate fixed point of the function $f$.

The algorithm that we use here, denoted by SFP-algorithm (\emph{simplicial
  fixed point algorithm}) for short, follows the path by using the
lexicographic pivoting rule from linear programming. The advantage is that
the so-called degeneration problem (i.e., the path ends up in a circuit)
cannot occur. We don't want to enlarge on that fact and will only give the
definition of lexicographically positive matrices. For a detailed
description we refer to \cite[Chapters~2\&3]{Van96}.

\begin{definition} A row vector is called \emph{lexicographically positive}, if its first nonzero entry is positive. A matrix $W$ is called \emph{lexicographically positive} denoted by $W \succ 0$, if every row vector is lexicographically positive. \end{definition}

\subsection{The SFP-algorithm}\label{subsec:Algo}

To compute a fixed point of a continuous function $\phi:\R^N\rightarrow
\R^N$ the SFP-algorithm uses a suitable homotopy
$\vartheta:\R^N\times[0,1] \rightarrow \R^N$ and a pivoting method to get
from an $N$-simplex $\tau^0 \subset \R^N \times\{0\}$ to an $N$-simplex
$\tau^1 \subset \R^N \times \{1\}$ which yields an approximate fixed point
of $\phi$. For this purpose we have to triangulate the set $\R^N \times
[0,1]$ suitably.

Let $\T$ be a triangulation of $\R^N \times [0,1]$ with the restriction
$\T^{0} \subset \R^N \times \{0,1\}$, i.e., $\T$ is two-layered. We denote
elements of $\R^N \times [0,1]$ by $ y = (v_{1}, \ldots ,v_{N},t)$ with $
v \in \R^N$ and $t \in [0,1]$ and define the \emph{projection onto the
  first factor} $p_1: \R^N \times [0,1] \rightarrow \R^N$, $p_1(v,t)=v.$
Suppose that the $N$-simplex $\tau=\langle y^1, y^2, \ldots,
y^{N+1}\rangle \in \T^N$. We define the diameter of the projection of
$\tau$ by 
\[\diam_{p}(\tau) := \max \{ \| p_1(y^i)-p_1(y^j)\| \ : \ i,j
\in \{1, \ldots, N+1\} \}.\]
Moreover, the mesh size of the projection of
$\T$ is defined by $$\mesh_{p}(\T):= \sup\{ \diam_{p}(\tau) \ : \ \tau \in
\T^N \}.$$ If $\tau=\langle y^1, \ldots, y^{N+1} \rangle \in \T^N$ and
$\tau \subset \R^N \times \{i\}$, $i\in\{0,1\}$, then $\tau_{p}:=\langle
p_1(y^1), \ldots, p_1(y^{N+1})\rangle$ is an $N$-simplex in $\R^N$. The
collection of all such simplices $\tau_{p}$ is denoted by $\T_{i}$.

We choose an arbitrary point $(c,0) \in \R^N \times [0,1]$ such that $(c,0)$ lies in the interior of an $N$-simplex $\tau^0 \in \T_{0}$. Consider the following homotopy mapping $\vartheta: \R^N \times [0,1] \rightarrow \R^N$ defined by 
\begin{equation}\label{eq:Merrill-homotopy-mapping}
\vartheta(v,t) = (1-t)c + t \phi(v).\end{equation}
A point $y$ is called a \emph{fixed point} of $\vartheta$, if $p_1(y) = \vartheta(y)$. Clearly, $(c,0)$ is the only fixed point of $\vartheta$ in $\R^N \times \{0\}$ and any fixed point $y$ of $\vartheta$ in $\R^N \times \{1\}$ projects to a fixed point of $\phi$, i.e., $p_1(y) = \phi(p_1(y))$. 
The concept of labelings establishes a way of studying the relation of the triangulation with approximate fixed points of $\phi$.

\begin{definition}\label{Def:Merril-LabelFunc} 
Let $\T$ be a two-layered triangulation of $\R^N \times [0,1]$. Then we define the \emph{labeling rule} $l: \R^N \times [0,1] \rightarrow \R^N$ by
\begin{equation}\label{eq:Merrill-Bezeichnungsfunktion}l(y)=\vartheta(y)-p_1(y).\end{equation}
Let the $N$-simplex $\tau=\langle y^1, \ldots, y^{N+1} \rangle \subset
\T^N$ be given. Then we call the $(N+1)\times(N+1)$ matrix
\begin{equation}\label{eq:Bezeichnungsmatrix} L(\tau):= \left(\begin{array}{ccc}1 & \hdots & 1 \\l(y^1) & \hdots & l(y^{N+1})\end{array}\right)\end{equation} 
the \emph{labeling matrix} of $\tau$.\end{definition}

The $N$-simplex $\tau$ is called \emph{complete}, if the system
\begin{equation}\label{eq:clsol} L(\tau) W = I_{N+1}, \quad W \succ 0\end{equation}
has a solution $W^* \in \R^{(N+1) \times (N+1)}$. Complete simplices play
an important role in the following since a complete $N$-simplex $\tau \subset
\R^N\times \{1\}$ contains an approximate fixed point of $\phi$. In
addition, by choosing the mesh size of the triangulation small enough we
can claim any accuracy of the approximate fixed point.

\begin{proposition}\label{satz:Merrill-ApproxFP} 
    Let $D$ be compact and $\phi: D \subset \R^N \rightarrow \R^N$ be
    continuous. For $\varepsilon>0$ let $\delta>0$ be such that for all
    $x,y \in D$ we have the implication $\|x-y\|<\delta \Rightarrow
    \|\phi(x)-\phi(y)\|<\varepsilon$. Let $\T$ be a two-layered
    triangulation of $\R^N \times [0,1]$ with $\mesh(\T)<\delta$ and
    $\tau=\langle y^1, \ldots, y^{N+1}\rangle \subset \R^N \times \{1\}$ a
    complete simplex in $\T$ with $y^j = (v^j,t_{j})$ for all $j=1,
    \ldots, N+1$. Let $\lambda \in \R^N_+$ be the solution of the system
\begin{equation}\label{eq:wstar}L(\tau)\tilde \lambda=e_{1}, \quad \tilde \lambda \in \R^{N+1}_{+}.\end{equation}
Then $v^* := \sum_{j=1}^{N+1} \lambda_j v^j$ is an approximate fixed point
of $\phi$, i.e., $\|\phi(v^*)-v^*\|<\varepsilon.$\end{proposition}

\begin{proof}
    Since $\tau \subset \R^N \times \{1\}$ we have $t_j=1$ for all $j=1,
    \ldots, N+1$ and so $ l(y^j) = \vartheta(y^j)-p_1(y^j) =
    \phi(v^j)-v^j$ by (\ref{eq:Merrill-homotopy-mapping}) and
    (\ref{eq:Merrill-Bezeichnungsfunktion}). Thus (\ref{eq:wstar}) is
    equivalent to
\begin{equation}\label{eq:sumlambda}
 (i) \quad \sum_{j=1}^{N+1} \lambda_j = 1 \qquad \text{and}\qquad (ii) \quad \sum_{j=1}^{N+1} \lambda_j \phi(v^j) = \sum_{j=1}^{N+1} \lambda_j v^j.
\end{equation}
By (\ref{eq:sumlambda})(i) $v^*$ is a convex combination of the $v^1, \ldots, v^{N+1}$, i.e., $v^* \in \tau$. But then we have $\|v^*-v^j\|<\delta$ for all $j=1,\ldots, N+1$ and by continuity of $\phi$ we have $\|\phi(v^*)-\phi(v^j)\|<\varepsilon$ for all $j=1,\ldots, N+1$. Together this yields
\begin{eqnarray*}
 \|\phi(v^*)-v^*\| &\stackrel{(\ref{eq:sumlambda})(i)}{=}& \| (\sum_{j=1}^{N+1} \lambda_j) \phi(v^*)-\sum_{j=1}^{N+1} \lambda_j v^j\| \\
	&\stackrel{(\ref{eq:sumlambda})(ii)}{=}& \| \sum_{j=1}^{N+1} \lambda_j \phi(v^*)-\sum_{j=1}^{N+1} \lambda_j \phi(v^j)\| \\
	&= & \sum_{j=1}^{N+1} \lambda_j \| \phi(v^*)- \phi(v^j)\| \stackrel{(\ref{eq:sumlambda})(i)}{<} \varepsilon.
\end{eqnarray*}
\end{proof}

To obtain a complete simplex in $\R^N \times \{1\}$ we first characterize
the complete simplices. To this end we define the graph $G(V,E)$ of all
complete simplices as follows. An $(N+1)$-simplex $\eta$ of $\T$ is a
\emph{node}, if it has at least one complete facet $\tau$. Two nodes are
\emph{adjacent} and connected by an edge, if they share a common complete
facet. The \emph{degree of a node} $\eta$ is the number of nodes
adjacent to $\eta$, denoted by $\deg(\eta)$.

Recall that $\tau^0$ is the $N$-simplex lying on $\R^N \times \{0\}$ and
containing $(c,0)$ in its interior. Let $\eta^0$ be the unique
$(N+1)$-simplex of $\T$ having $\tau^0$ as its facet. Then we have
(cf. \cite[Lemma 4.3.3, \ Lemma 4.3.4 \ and Theorem 4.3.5]{Yang}).

\begin{lemma}\label{lem:tauvollst} The $N$-simplex $\tau^0$ is the only
    complete simplex on $\R^N \times \{0\}$. \end{lemma}

\begin{lemma}\label{lem:degreetau} Given the graph $G(V,E)$ defined as above, for each node $\eta$ of $G(V,E)$, we have
\begin{enumerate}
  \item if $\eta$ has a complete facet lying on $\R^N \times \{0\}$ or
    $\R^N \times \{1\}$, then $\deg(\eta)=1$;
\item in all other cases, $\deg(\eta)=2$.
\end{enumerate}
\end{lemma}

\begin{theorem}\label{theo:PathPossibility} 
For the graph  $G(V,E)$ defined as above, each connected component of $G(V,E)$  has one of the following five forms\begin{enumerate}
\item a simple circuit (i.e. a path $(e_{0,1}, e_{1,2}, \ldots e_{k-1,k})$, $k\in \N$ with $e_{0,1}=e_{k-1,k}$ and $e_{i,i+1}\neq e_{j,j+1}$ for $i\neq j$ and $i,j \in \{1, \ldots, k-1\}$);
\item a finite simple path (i.e. a path without circuits) whose two end
  nodes all have a complete facet lying on $\R^N \times \{1\}$;
\item an infinite simple path starting with an $(N+1)$-simplex which has a
  complete facet lying on $\R^N \times \{1\}$;
\item a finite simple path which starts with the $(N+1)$-simplex
  $\eta^0$ and ends with another $(N+1)$-simplex having a complete facet
  on $\R^N \times \{1\}$;
\item an infinite simple path starting with the $(N+1)$-simplex $\eta^0$.
\end{enumerate}
\end{theorem}

From the point of view of computation we are interested in case $(iv)$. In
this case we can algorithmically go from $\eta^0$ to a simplex $\eta^*
\in \R^N \times \{1\}$ containing an approximate fixed point of $\phi$ by
Proposition \ref{satz:Merrill-ApproxFP}.  A schematic description is given
in Figure \ref{fig:AlgoMerrill}.

\begin{minipage}{\textwidth}
	\singlespacing 
	
	\textsc{The Simplicial Fixed Point Algorithm} 
	\vspace{0.3em}
	\hrule \hrule \vspace{0.8em}
	\setlength{\leftskip}{2em} \setlength{\rightskip}{2em}

	\textbf{Step (0) } Set $\mesh_{p}(\T)< \delta$. Let $\tau^0$ be the unique $N$-simplex of $\T^N$ containing $(c,0)$ in its interior. Let $\eta^0$ be the unique $(N+1)$-simplex in $\T$ which has $\tau^0$ as its facet. Let $y^+$ be the vertex of $\eta^0$ that is not a vertex of $\tau^0$. Set $k=0$.

\medskip
	\textbf{Step (1) } Compute $W=L^{-1}(\tau^k)$ with $L$ from $(\ref{eq:Bezeichnungsmatrix})$. Let $W_{i}$ denote the $i^{th}$ row of $W$. Compute $l(y^+)$ and let $q=(1, l^\top(y^+))^\top$. Let $p=(p_{1}, \ldots, p_{n+1})^\top=Wq$ denote the coefficient vector of the linear combination $l(y^+)= \sum_{i=1}^{N+1} p_{i}l(y^i)$. Compute $\zeta \in \{1,\ldots,N+1\}$ so that the quotient
	$$\frac{W_{\zeta}}{p_{\zeta}}= \min_{\prec} \left \{ \frac{W_{h}}{p_{h}} \ : \ p_{h}>0 , \ h =1,\ldots, N+1\right\}$$
	is lexicographically positive minimal. Note that $\zeta$ is unique
        (cf. \cite[Theorem 4.2.7]{Yang}). Let $\tau^{k+1}$ be the facet of
        $\eta^k$ opposite $y^\zeta$. If $\tau^{k+1}$ lies on $\R^N
        \times \{1\}$, this facet yields an approximate fixed point of
        $\phi$ and stop. If $\tau^{k+1}$ does not lie on $\R^N \times
        \{1\}$, go to Step (2).

\medskip
	\textbf{Step (2) } Find a simplex $\eta^{k+1}$ sharing the facet $\tau^{k+1}$ with $\eta^k$ (which is unique by Lemma \ref{lem:degreetau}), and let $y^+$ be the vertex $\eta^{k+1}$ not being a vertex of $\tau^{k+1}$. Set $k=k+1$ and return to Step (1).\\	
	\hrule \hrule
	\setlength{\leftskip}{0em} \setlength{\rightskip}{0em}
	\onehalfspacing 
\end{minipage}
\begin{figure}[h]
\caption{The SFP-algorithm} 
\label{fig:AlgoMerrill}
\end{figure} 

\begin{remark}\label{rem:Merrill}
    In order to guarantee case $(iv)$ in Theorem
    $\ref{theo:PathPossibility}$ Merrill (cf. \cite{Yang}) gave a
    contraction condition that is sufficient for the convergence of the
    algorithm for a u.s.c. point-to-set mapping $\phi: \R^N \rightarrow
    \R^N$, in particular for a continuous single-valued function $\phi:
    \R^N \rightarrow \R^N$, see \cite[Theorem 4.3.6]{Yang}. The method of
    proof is to show that there is only a compact subset $D \subset \R^N$
    yielding complete simplices, so the path must be finite. Note that
    this condition is not sufficient for convergence, if we choose
    $\phi: \R^N_+ \rightarrow \R^N$. In particular the function $\phi$
    defined in $(\ref{eq:phix})$ satisfies Merrill's condition
    (cf. \cite[Satz 4.28]{Geiselhart}) but we have to impose other
    conditions to guarantee convergence.
\end{remark}

\subsection{Using the SFP-algorithm for computing decay points}\label{subsec:Anpassung}

Now we want to use the SFP-algorithm to compute a decay point $w \in
\Omega(\Gammu)$ of the monotone operator $\Gammu: \R^N_{+} \rightarrow
\R^N_{+}$ which satisfies $\Gammu(0)=0$.

In the following the aim is to find a suitable function $\phi$ whose fixed
points $w = \phi(w)$ correspond to decay points $w \in \Omega(\Gammu)$,
and to show that the SFP-algorithm converges for this choice of
$\phi$. Since $\Gammu(0)=0$ and this point yields no information, we have
to exclude $0$ from being a fixed point of $\phi$. Also, in order to show
that complete simplices can only lie in a compact subset of $\R^N_+$ it is
desirable to have $\phi$ \emph{small} for \emph{large} $v$.

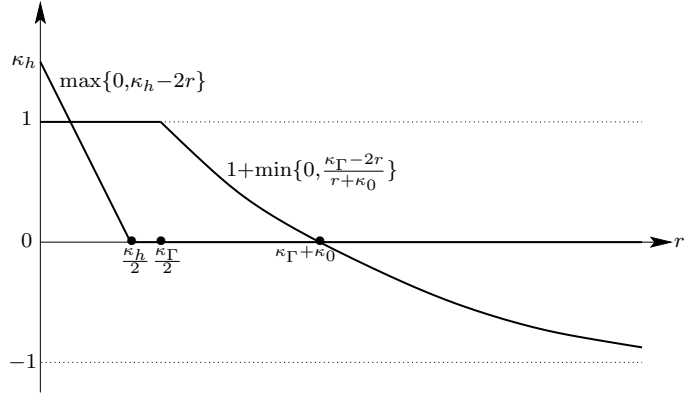
\begin{figure}[ht] 
\centering
\begin{overpic}[scale=0.4]{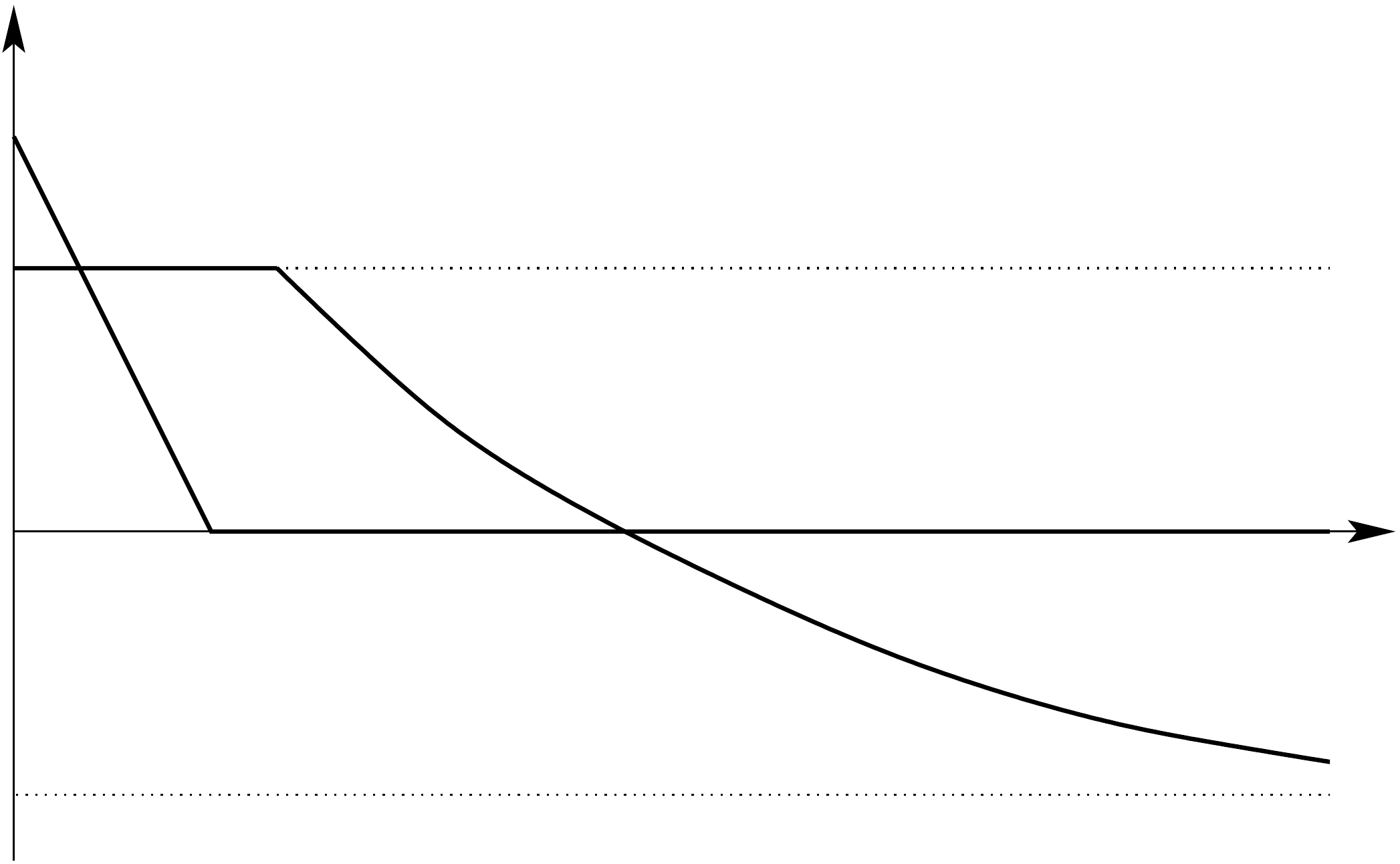}
\put(4,48){${\scriptstyle \max\{0,\kappa_{h}-2r\}}$} 
\put(30,34){${\scriptstyle 1+\min\{ 0, \frac{\kappa_{\Gamma}-2r}{r+\kappa_{0}}\}}$}
\put(-4,4){${\scriptstyle -1}$} 
\put(-2,23){${\scriptstyle 0}$} 
\put(-2,42){${\scriptstyle 1}$} 
\put(-3.5,51.5){${\scriptstyle \kappa_{h}}$}
\put(14.4,23.1){${\scriptstyle \bullet}$} 
\put(13.4,20.1){${\scriptscriptstyle \frac{\kappa_{h}}{2}}$}
\put(19,23.1){${\scriptstyle \bullet}$}  
\put(18.4,20.1){${\scriptscriptstyle \frac{\kappa_{\Gamma}}{2}}$}
\put(43.8,23.1){${\scriptstyle \bullet}$}  
\put(37.6,21.5){${\scriptscriptstyle \kappa_{\Gamma}+\kappa_{0}}$}
\put(100, 23){${\scriptstyle r}$}
\end{overpic}
\caption{The components of the function $\phi$}
\label{fig:KompPhi}
\end{figure}

Consider the function $\phi: \R^N_{+} \rightarrow \R^N$ defined by
\begin{equation}\label{eq:phix} \phi(v) = \Gammu(v)\left(1+ \min\left\{ 0, \frac{\kappa_{\Gamma}-2\|v\|}{\|v\|+\kappa_{0}}\right\} \right) + \max\left\{0,\kappa_{h}-2\|v\|\right\}e.\end{equation}
Here let $\kappa_{0}>0$, $\kappa_{\Gamma}>\kappa_{h}>0$ and
$e:=\sum_{i=1}^N e_i$ the $N$-dimensional vector of ones. We illustrate
the components of $\phi$ in Figure \ref{fig:KompPhi}.

Some properties of $\phi$ are as follows:
\begin{enumerate}
  \item $\phi$ is continuous on $\R^N_{+}$ since $\Gammu$ is continuous on
    $\R^N_{+}$.
\item For large $v \in \R^N_{+}$ it holds $\phi(v)<0$.
\item It holds $\phi(0)= \kappa_{h}e \gg 0$, i.e., the origin cannot be a
  fixed point of $\phi$.
\end{enumerate}

\begin{figure}[ht]
\centering
\begin{overpic}[scale=.35]{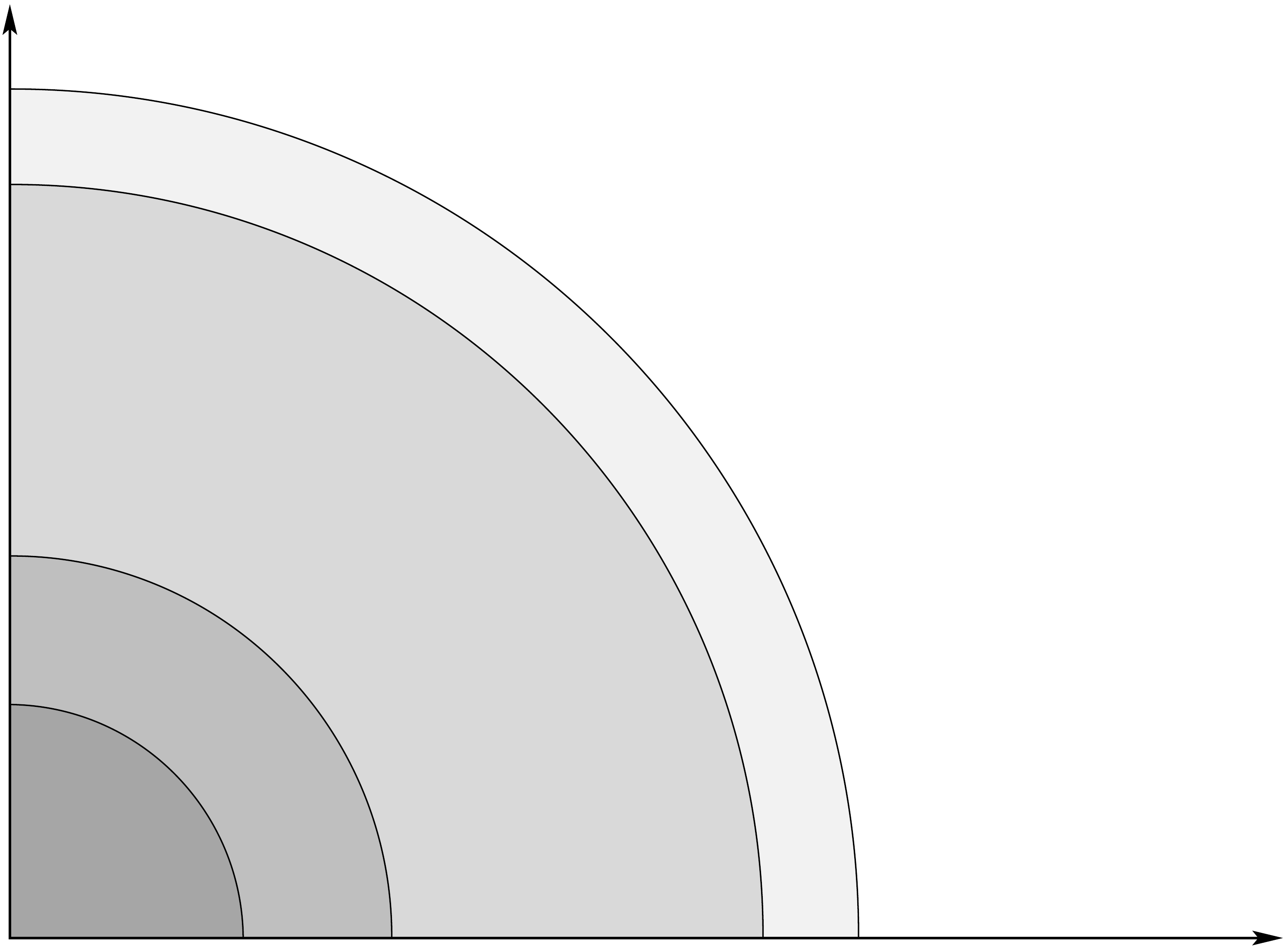} 
\put(8,6){${\scriptstyle \1}$} 
\put(20.5,11){${\scriptstyle \2}$}  
\put(39,17){${\scriptstyle \3}$} 
\put(57,22){${\scriptstyle \4}$} 
\put(66,24){${\scriptstyle \5}$}
\put(59.8,65){${\scriptstyle \1}$} 
\put(65,65){${\scriptstyle \phi(v)=\Gammu(v)+ (\kappa_{h}-2\|v\|)e \gg 0}$}
\put(59,60){${\scriptstyle \2}$} 
\put(65,60){${\scriptstyle \phi(v)=\Gammu(v)>0}$}
\put(58.5,55){${\scriptstyle \3}$} 
\put(65,55){${\scriptstyle \phi(v)=\Gammu(v)\left(1+ \frac{\kappa_{\Gamma}-2\|v\|}{\|v\|+\kappa_{0}}\right)>0}$}
\put(57,50){${\scriptstyle \4,\5}$} 
\put(65,50){${\scriptstyle \phi(v)=\Gammu(v)\left(1+ \frac{\kappa_{\Gamma}-2\|v\|}{\|v\|+\kappa_{0}}\right)<0}$}
\put(17.5,-2.5){${\scriptstyle \frac{\kappa_{h}}{2}}$}
\put(29,-2.5){${\scriptstyle \frac{\kappa_{\Gamma}}{2}}$}
\put(54.5,-1.5){${\scriptscriptstyle \kappa_{\Gamma}+\kappa_{0}}$}
\put(65.5,-1.5){${\scriptscriptstyle \kappa_{\Gamma}+\kappa_{0}+\delta}$}
\end{overpic}
\caption{The definition of $\phi$ illustrated on the positive orthant} 
\label{fig:posOrth}
\end{figure}
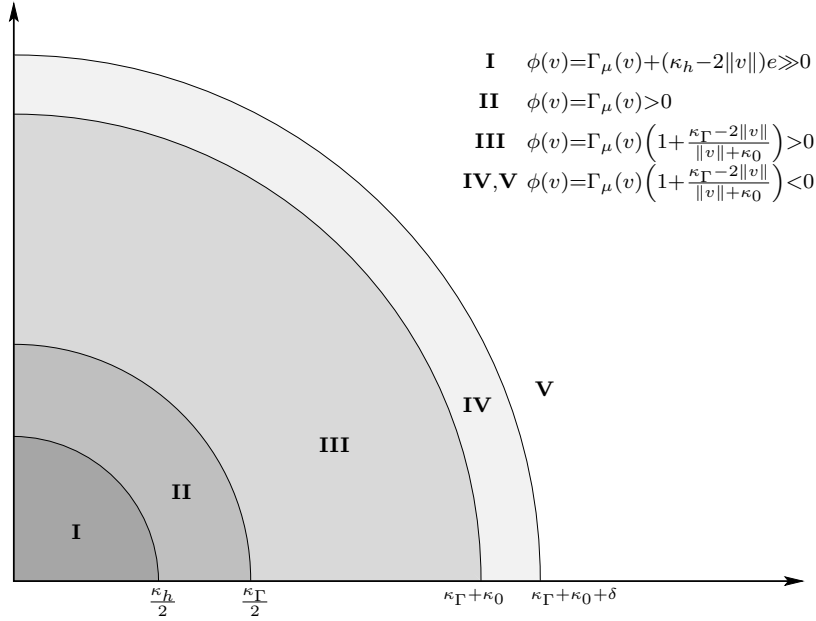

In Figure \ref{fig:posOrth} we illustrate the definition of $\phi$ on the positive orthant. To this end we partition the positive orthant in five regions:

$$\begin{array}{lcl}
\begin{array}{ccl}
\1 &=& \left\{ v \in \R^N_{+} :  \| v\|  \in [0, \kappa_{h}/2) \right\},\\
\2 &=& \left\{ v \in \R^N_{+} :  \| v\|  \in [\kappa_{h}/2, \kappa_{\Gamma}/2) \right\},\\
\3 &=& \left\{ v \in \R^N_{+} :  \| v\|  \in [\kappa_{\Gamma}/2, \kappa_{\Gamma}+\kappa_{0}) \right\},\\
\4 &=& \left\{ v \in \R^N_{+} :  \| v\|  \in [\kappa_{\Gamma}+\kappa_{0},\kappa_{\Gamma}+\kappa_{0}+\delta] \right\},\\
\5 &=& \left\{ v \in \R^N_{+} :  \| v\|  \in (\kappa_{\Gamma}+\kappa_{0}+\delta, \infty) \right\},
\end{array}
& &
\begin{array}{ccl}
\1' &=& \1 \times[0,1],\\
\2' &=& \2 \times[0,1],\\
\3' &=& \3 \times[0,1],\\
\4' &=& \4 \times[0,1],\\
\5' &=& \5 \times[0,1].
\end{array}\end{array}$$

\bigskip

\bigskip

The next proposition indicates the relation between fixed points of $\phi$
and decay points of $\Gammu$.

\begin{proposition}\label{satz:FPAP}
    Let $\phi: \R^N_+ \rightarrow \R^N$ be defined as in $(\ref{eq:phix})$
    and assume that $\Gammu: \R^N_{+} \rightarrow \R^N_{+}$ is monotone
    and satisfies the small gain condition $(\ref{eq:Merrill:SGB})$. Let
    $s \in \R^N_{+}$ be a fixed point of the function $\phi$, i.e., $s =
    \phi(s)$. Then $s$ lies in the set of decay of the function $\Gammu$,
    i.e., $s \in \Omega(\Gammu)$. Moreover, $s \in \1$.
\end{proposition}

\begin{proof}
We distinguish between the following two cases for $s \in \R^N_{+}$:
\begin{enumerate}
\item $0 \leq \|s\| < \frac{\kappa_{h}}{2}$:
	In this case $s=\phi(s)=\Gammu(s)+(\kappa_{h}-2\|s\|)e \gg \Gammu(s)$. It follows that $s \gg \Gammu(s)$, so $s \in \Omega(\Gammu)$. In particular, $s\in \1$.
\item $\frac{\kappa_{h}}{2} \leq \|s\|$:
	 We have $s=\phi(s)=\Gammu(s)(1+ \min\{0,\frac{\kappa_{\Gamma}-2\|s\|}{\|s\|+\kappa_{0}}\}) \leq \Gammu(s) $ but this is a contradiction to the small gain condition $(\ref{eq:Merrill:SGB})$, so this case cannot occur. \end{enumerate}
\end{proof}

In the following we will always use the $\tilde K_{1}(\delta)$-triangulation. This triangulation has the essential advantage that the vertices of an $N$-simplex $\tau=\langle y^1, \ldots, y^{N+1}\rangle$ are in the order of $\R^{N+1}_+$, i.e., it holds $ y^1< \ldots < y^{N+1}.$ Note that $y=(v,t) \in \R^{N}_{+} \times \{0,1\}$.\\
Again, the SFP-algorithm starts with the $(N+1)$-simplex $\eta^0$ which
has the $N$-simplex $\tau^0 \in \R^N \times \{0\}$ as a facet containing
$(c,0)$ in its interior, where $c$ determines the homotopy mapping $\vartheta$ in \eqref{eq:Merrill-homotopy-mapping}.  
Here we choose $c \in \1 \cup \2$ and any
approximate fixed point $c'$ will also lie in $\1 \cup \2$, see Theorem
\ref{satz:nichtnachIV}. Then the algorithm follows the path of complete
$N$-simplices.  If we can show that this path is finite and inside of the
positive orthant, then we get, by Theorem \ref{theo:PathPossibility}, that
the SFP-algorithm ends up with a $(N+1)$-simplex containing a complete
facet on $\R^N_{+} \times \{1\}$.
Proposition \ref{satz:Merrill-ApproxFP} now tells us that this simplex contains an approximate fixed point of $\phi$.\\

A first rough estimation where the path of complete simplices can run is given in the next proposition.

\begin{proposition}\label{satz:nichtnachV}
Let $\phi: \R^N_+ \rightarrow \R^N$ be defined as in $(\ref{eq:phix})$ and assume that $\Gammu: \R^N_{+} \rightarrow  \R^N_{+}$ is monotone. Assume that the constant $c$ used in \eqref{eq:Merrill-homotopy-mapping} satisfies $c\in \1\cup\2$ and let $\tau=\langle y^1,\ldots, y^{N+1} \rangle$ be an $N$-simplex in $\5'$. Then $\tau$ is not complete.
\end{proposition}

\begin{proof}
 We prove this by contradiction. Assume that $\tau$ is complete.  Then the linear system 
\begin{equation}\label{eq:linsysnichtV}
 L(\tau)W=I_{N+1}, \quad W \succ 0
\end{equation}
with $L$ defined as in (\ref{eq:Bezeichnungsmatrix}), has a lexicographically positive solution $W$. We have $y^j = (v^j,t_j) \in \R^N_+ \times \{0,1\}$ and by $\tau \subset \5'$ we have $\|v^j\|> \kappa_\Gamma+\kappa_0+\delta$, i.e., $\phi(v^j)<0$ for all $j =1, \ldots, N+1$. So we have the following two cases for $$l(y^j) = \vartheta(v^j,t_j) - v^j = (1-t_j)c+ t_i\phi(v^j)-v^j.$$
\begin{enumerate}
 \item If $t_j=1$, then $l(y^j)\ll 0$;
 \item If $t_j=0$, then $l(y^j) = c-v^j$. Since $c \in \1 \cup \2$ and $v^j \in \5$ for all $j=1,\ldots, N+1$ there exists a component $i^* \in \{1, \ldots, N\}$ with $c_{i^*} < v^1_{i^*} \leq v^j_{i^*}$ for all $j=2,\ldots, N+1$.
\end{enumerate}
Together it follows 
\begin{equation}\label{eq:LabelnichtV}
 l(y^j)_{i^*}<0 \qquad \text{for all } j=1, \ldots, N+1.
\end{equation}
Let $L_l$ denote the $l^{th}$ row of $L$ and let $W^m$ denote the $m^{th}$ column of $W$. Then we have $W^1 \in \R^N_+ \backslash \{0\}$ since $W$ is lexicographically positive. By (\ref{eq:LabelnichtV}) we have $L_{i^*+1}\ll 0$. But then $L_{i^*+1}W^1<0$ in contradiction to (\ref{eq:linsysnichtV}). So $\tau$ cannot be complete.
\end{proof}

Note that this does not show that the path starting in $\eta^0$ is
inside the positive orthant. To prove this we have to look at the boundary
of the positive orthant. Here we need some additional assumptions. Note
that for $\Gammu$, $\Gamma \in (\Ki \cup \{0\})^{N \times N}$ is the
underlying gain matrix (cf. Section~\ref{subsec:Grundlagen-Matrizen}), and
by Remark~\ref{remark:ground-Gammu} the operator $\Gammu$ is monotone.

\begin{theorem}\label{satz:intOrth} 
Let $\phi: \R^N_+ \rightarrow \R^N$ be defined as in $(\ref{eq:phix})$ and assume that the underlying gain matrix $\Gamma$ is irreducible. Let $\tau=\langle y^1,\ldots, y^{N+1} \rangle$ be an $N$-simplex on the boundary of the positive orthant. If $\|v^{N+1}\| = \|p_1(y^{N+1})\|<\kappa_{0}+\kappa_{\Gamma}$ then $\tau$ is not complete. \end{theorem}

\begin{proof}
If $\tau$ is an $N$-simplex on the boundary of the positive orthant then there exists an index $i^*\in \{1, \ldots, N\}$ with $v^j_{i^*}=0$ for all $j=1, \ldots, N+1$. We prove by contradiction that $\tau$ cannot be complete, if $\|v^{N+1}\| <\kappa_{0}+\kappa_{\Gamma}$. So assume 
\begin{equation}\label{eq:linsysintOrth}
 L(\tau)W=I_{N+1},\qquad W \succ 0
\end{equation}
has the solution $W^*$ and let $\lambda \in \R^{N+1}_{+}$ be the first column of $W^*$. Then it follows by $(\ref{eq:linsysintOrth})$ and using $(\ref{eq:Bezeichnungsmatrix})$ that
\begin{equation}\label{eqn:(1)} 
 \sum_{j=1}^{N+1}\lambda_{j}l(y^j) = 0, \quad \sum_{j=1}^{N+1}\lambda_{j} = 1. 
\end{equation}
The case $\|v^{N+1}\|<\frac{\kappa_{h}}{2}$ yields, using $(\ref{eq:Merrill-Bezeichnungsfunktion})$, $l(y^j)_{i^*}=(1-t_{j})c_{i^*}+t_{j}\phi(v^j)_{i^*}>0$ for all $j\in \{1, \ldots, N+1\}$, so $\sum_{j=1}^{N+1}\lambda_{j}l(y^j)_{i^*}>0$ since $\lambda\neq 0$ and $c \gg 0$ (since $c$ lies in the interior of a simplex $\tau$). But this is a contradiction to $(\ref{eqn:(1)})$.\\
Now assume $\|v^{N+1}\|<\frac{\kappa_{\Gamma}}{2}$. Then it holds $l(y^j)_{i^*}\geq0$. In particular, $l(y^j)_{i^*}=(1-t_{j})c_{i^*} + t_{j}\Gammu(v^j)_{i^*}+ t_{j}(\kappa_{h}-2\|v^j\|)>0$ for $\|v^j\|<\frac{\kappa_{h}}{2}$. 
Let $r \in \{1, \ldots, N+1\}$ with $t_{1}= \ldots = t_{r}=0$ and $t_{r+1}= \ldots = t_{N+1}=1$ as well as 
$\rho \in \{1, \ldots, N+1\}$ with $\|v^\rho\| <\frac{\kappa_{h}}{2}$ and $\|v^{\rho+1}\| \geq \frac{\kappa_{h}}{2}$. Then equation $(\ref{eqn:(1)})$ is equivalent to
	\begin{equation}\label{eqn:(2)} \sum_{j=1}^{r}\lambda_{j}c + \sum_{j=r+1}^{N+1}\lambda_{j} \Gammu(v^j)+ \sum_{j=r+1}^\rho\lambda_{j}(\kappa_{h}-2\|v^j\|)e = \sum_{j=1}^{N+1}\lambda_{j}v^j, \quad \sum_{j=1}^{N+1}\lambda_{j} = 1. \end{equation}
Since $v^j_{i^*}=0$ for all $j=1, \ldots, N+1$ it follows $\lambda_{j}=0$ for $j=1, \ldots, \hat r$ with $\hat r:= \max \{ r, \rho\}$ and equation (\ref{eqn:(2)}) is equivalent to
	\begin{equation}\label{eqn:(3)} \sum_{j=\hat r+1}^{N+1}\lambda_{j}\Gammu(v^j) = \sum_{j=\hat r+1}^{N+1}\lambda_{j}v^j, \quad \sum_{j=\hat r+1}^{N+1}\lambda_{j} = 1. \end{equation}
Now there exists a largest index $\tilde r \in \{1, \ldots, N+1\}$ with $\Gammu(v^{\tilde r})_{i^*}=0$ and $\Gammu(v^{\tilde r})_{i^*}>0$. Since the $v^j$ are ordered by the $\tilde K_{1}(\delta)$-triangulation, it follows by monotonicity of $\Gammu$,
\begin{equation}\label{eq:GammumonotonintOrt}
 \Gammu(v^1)< \ldots< \Gammu(v^{\tilde r})< \ldots< \Gammu(v^{N+1}),
\end{equation}
and thus $\lambda_{\tilde r+1}= \ldots = \lambda_{N+1}=0$ which leads to
	\begin{equation}\label{eqn:(4)} \sum_{j=\hat r+1}^{\tilde r}\lambda_{j}\Gammu(v^j) = \sum_{j=	\hat r+1}^{\tilde r}\lambda_{j}v^j, \quad \sum_{j=\hat r+1}^{\tilde r}\lambda_{j} = 1.\end{equation}
Without loss of generality assume $v^{\tilde r}=[v^{\tilde r}_{1}, \ldots, v^{\tilde r}_{l}, 0, \ldots, 0]^\top$ with $v^{\tilde r}_{j}>0$ for $j=1, \ldots, l$, $l \leq N$. Equation $(\ref{eqn:(4)})$ implies $\Gammu(v^{\tilde r}) = [*_{1}, \ldots, *_{l}, 0, \ldots, 0]^\top$ with $*_{j}\geq0$, $j=1, \ldots, l$. But then $\Gamma$ is of the form 
	\begin{equation}\label{eqn:(5)} \Gamma=\left(\begin{array}{cc}\Gamma_{11} & \Gamma_{12} \\0 & \Gamma_{22}\end{array}\right)\end{equation} 
	with $\Gamma_{11} \in \R_{+}^{l\times l}$ and $\Gamma_{22} \in \R_{+}^{(N+1-l)\times (N+1-l)}$. This means that $\Gamma$ is reducible, a contradiction to the assumption. So this case cannot occur.\\
Now assume $\|v^{N+1}\|<\kappa_{0}+\kappa_{\Gamma}$. Define $\tilde\Gamma(v^j):=\Gammu(v^j)\left(1+\frac{\kappa_{\Gamma}-2\|v^j\|}{\|v^j\|+\kappa_{0}}\right)$ if $\|v^j\|>\frac{\kappa_{\Gamma}}{2}$. For $\|v^j\|<\kappa_{0}+\kappa_{\Gamma}$ it holds $\tilde \Gamma(v^j)_{k}\geq0$ and $\Gammu(v^j)_{k}=0 $, if and only if $\tilde \Gamma (v^j)_{k}=0$. 
The same argumentation as above  provides
	\begin{equation} %\tag{\ref{eqn:(4)}$^*$}\label{eqn:(4*)} 
\sum_{j=\hat r+1}^{\tilde r}\lambda_{j}\tilde \Gamma(v^j) = \sum_{j=	\hat r+1}^{\tilde r}\lambda_{j}v^j, \quad \sum_{j=\hat r+1}^{\tilde r} \lambda_{j} = 1.\end{equation}
With $v^{\tilde r}=[v^{\tilde r}_{1}, \ldots, v^{\tilde r}_{l}, 0, \ldots, 0]^\top$, $v^{\tilde r}_{j}>0$ for $j=1, \ldots, l$ it follows $\tilde \Gamma(v^{\tilde r}) = [\tilde*_{1}, \ldots, \tilde*_{l}, 0, \ldots, 0]^\top$ with $\tilde*_{j}\geq 0$, $j=1, \ldots, l$. All in all we get $\tilde \Gamma(v^{\tilde r})<\Gammu(v^{\tilde r})= [*_{1}, \ldots, *_{l}, 0, \ldots, 0]^\top$ with $*_{j}>\tilde*_{j}\geq0$, $j=1, \ldots, l$. But then $\Gamma$ is of the form $(\ref{eqn:(5)})$, a contradiction to the assumption.
\end{proof}

\begin{figure}[ht] 
\centering
\begin{overpic}[scale=.35]{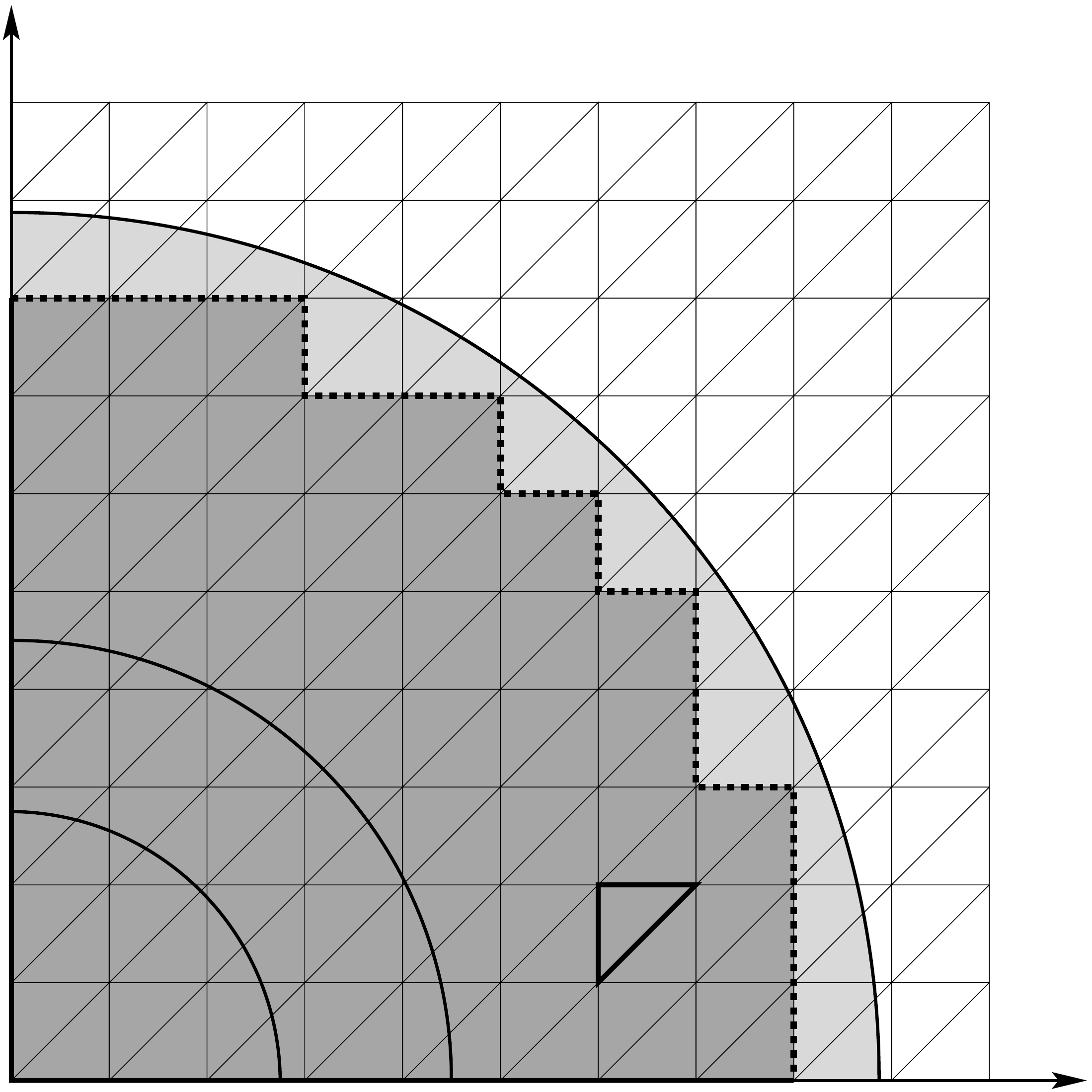} 
\put(13,13){${\scriptstyle \1'}$} 
\put(22,22){${\scriptstyle \2'}$} 
\put(38,40){${\scriptstyle \3'}$} 
\put(65,65){${\scriptstyle \4'}$} 
\put(56,15){${\scriptstyle \eta}$}
\put(64.2,32){${\scriptstyle \tau}$}
\put(23,-2.5){${\scriptscriptstyle \frac{\kappa_{h}}{2}}$}
\put(39,-2.5){${\scriptscriptstyle \frac{\kappa_{\Gamma}}{2}}$}
\put(76,-1.5){${\scriptscriptstyle \kappa_{\Gamma}+\kappa_{0}}$}
\end{overpic}
\caption{$\tilde K_{1}(\delta)$-triangulation and the maximum region of the path ($\tau$-simplices are 1-dimensional and $\eta$-simplices are 2-dimensional)}
\label{fig:posOrthTriang}
\end{figure}
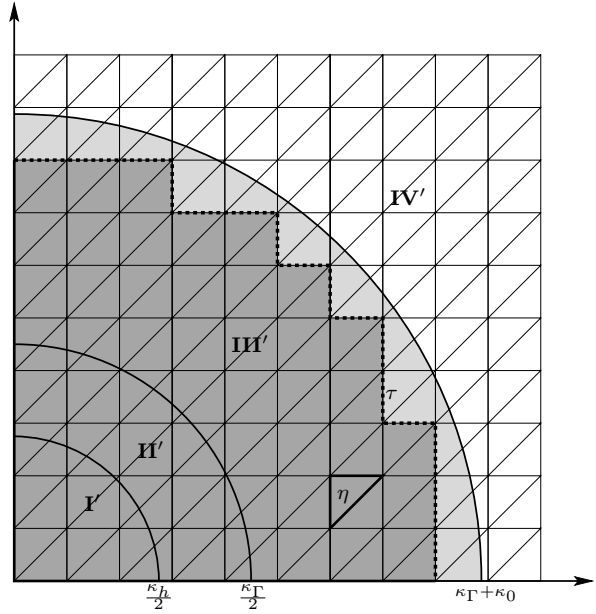

In other words Theorem \ref{satz:intOrth} provides that no $N$-simplex
$\tau \in \1' \cup \2' \cup \3'$ lying on the boundary of the positive
orthant can be complete. So it remains to show that the path starting in
$\eta^0$ cannot enter the set $\4'$. For this purpose we show in the next
theorem that the path of complete simplices runs inside of the region
which is painted dark grey in Figure \ref{fig:posOrthTriang}. To prove
this we demand an upper bound for the feasible size of $\delta$.

\begin{theorem}\label{satz:nichtnachIV} 
    Let $\phi: \R^N_+ \rightarrow \R^N$ be defined as in $(\ref{eq:phix})$
    and assume that $\Gammu: \R^N_{+} \rightarrow \R^N_{+}$ satisfies the
    small gain condition $(\ref{eq:Merrill:SGB})$. Furthermore assume that
    the underlying gain matrix $\Gamma$ does not contain any zero row,
    i.e., $\Gammu(e)_{i}\neq 0$ for all $i=1, \ldots, N$. Then there
    exists a $\delta>0$ such that for all simplices $\tau=\langle y^1,
    \ldots, y^{N+1} \rangle \subset \3'$ with $y^{N+2}:= y^1 + [\delta,
    \ldots, \delta,1]^\top \in \4'$ it follows that $\tau$ cannot be
    complete.

    In particular, any approximate solution $c'$ satisfies $c' \in \1 \cup
    \2$.
\end{theorem}

\begin{proof} Simplices satisfying $\tau\subset \3'$ with $y^{N+2}:= y^1 +
    [\delta, \ldots, \delta,1]^\top \in \4'$ are marked as black dotted
    lines in Figure \ref{fig:posOrthTriang}. We show that such a simplex
    $\tau$ cannot be complete, i.e., the system
\begin{equation}\label{eq:vollstIII} L(\tau)W=I_{N+1}, \quad W \succ 0\end{equation}
has no solution. %\\
First it holds for all $s \in \R^N_{+}$ with $\|s\| < \kappa_{\Gamma}+\kappa_{0}$ that
\begin{eqnarray*} \Gammu(s) &<& \max \left\{ \Gammu(v)  : v \in \R^N_{+} \wedge \|v\| = \kappa_{\Gamma}+\kappa_{0} \right\} \\
&\ll& \Gammu( [\kappa_{\Gamma}+\kappa_{0}, \ldots, \kappa_{\Gamma}+\kappa_{0}]^\top) =: \Gammu^{\max}. 
\end{eqnarray*}
Choose $k \in \N$ such that
\begin{equation}\label{eq:Wahlvonk} \frac{1}{2k-1} \Gammu^{\max} < c \end{equation} 
and set $\delta > 0$ such that
\begin{equation}\label{eq:Wahlvondelta} \delta < \min \left\{ \frac{\kappa_{\Gamma}-\kappa_{h}}{2\sqrt{N}} , \frac{\kappa_{\Gamma}+2\kappa_{0}}{2k\sqrt{N}} \right\}. \end{equation}
Now any approximate solution $c'$ can only be in $\1 \cup \2$ and not in
$\3$ since $\delta < \frac{\kappa_{\Gamma}-\kappa_{h}}{2\sqrt{N}}$. This
follows by Proposition \ref{satz:FPAP} because any fixed point $s^*$ of
$\phi$ has the property $\|s^*\| < \frac{\kappa_{h}}{2}$ and then
$$\|c'\|<\|s^*\|+\mesh_p(\tilde K_1(\delta)) = \|s^*\|+\delta \sqrt{N} \stackrel{(\ref{eq:Wahlvondelta})}{<} \frac{\kappa_{h}}{2} + \frac{\kappa_{\Gamma}-\kappa_{h}}{2\sqrt{N}}\sqrt{N} = \frac{\kappa_{\Gamma}}{2}.$$
The idea now is the following:
Let $r \in \{1, \ldots, N\}$ with $t_{r}=0$ and $t_{r+1}=1$. Then it holds
$$L(\tau) = \left(\begin{array}{cccccc}1 & \hdots & 1 & 1 & \hdots & 1 \\c-v^1 & \hdots & c-v^r & \phi(v^{r+1})-v^{r+1} & \hdots & \phi(v^{N+1})-v^{N+1}\end{array}\right).$$
Moreover, $v^1 \in \3$. So there exists at least one index $l \in \{1, \ldots, N+1\}$ with $(c-v^1)_{l}<0$. Then
\begin{equation}\label{eq:proof:c-vj} (c-v^j)_{l}<0 \quad \text{ for } \ j =1, \ldots, r. \end{equation}
Now $\phi(v)$ converges to zero, if $v$ tends to the boundary of $\4$, i.e., if $v \rightarrow \tilde v$ with $\|\tilde v\|=\kappa_{\Gamma}+\kappa_{0}$. So the aim is to get $\phi(v)$ as small as $(\phi(v^j)-v^j)_{l}<0$ holds for all $j=r+1, \ldots, N+1$.\\
Note that the function
$$g: \R^+ \rightarrow \R,\quad  a \mapsto 1 + \frac{\kappa_{\Gamma}-2a}{a+\kappa_{0}}$$
is strictly decreasing for $\kappa_{0}>0, \kappa_{\Gamma}>0$.
From the relation $v^{N+2} = v^1 + [\delta, \ldots \delta]^\top$ we get
$$ \kappa_{\Gamma}+\kappa_{0}< \| v^{N+2}\| = \| v^1 + [\delta, \ldots, \delta]^\top \| \leq \|v^1\| + \delta \sqrt{N}.$$
Under this assumption it follows for all $a > \kappa_{\Gamma}+\kappa_{0}-\delta \sqrt{N}$
\begin{eqnarray*}
g(a) &<& 1 + \frac{\kappa_{\Gamma}- 2 (\kappa_{\Gamma}+\kappa_{0}-\delta \sqrt{N})}{\kappa_{\Gamma}+\kappa_{0}-\delta\sqrt{N} + \kappa_{0}} \\
&=& \frac{\kappa_{\Gamma} + 2 \kappa_{0} - \delta\sqrt{N} + \kappa_{\Gamma} - 2 \kappa_{\Gamma} - 2 \kappa_{0}+ 2 \delta \sqrt{N}}{ \kappa_{\Gamma} + 2 \kappa_{0}- \delta \sqrt{N}} \\
&=& \frac{\delta \sqrt{N}}{\kappa_{\Gamma} + 2 \kappa_{0}- \delta \sqrt{N}} \stackrel{(\ref{eq:Wahlvondelta})}{<} \frac{\delta\sqrt{N}}{2k\delta\sqrt{N}-\delta\sqrt{N}} = \frac{1}{2k-1}.
\end{eqnarray*}
Together with $(\ref{eq:Wahlvonk})$ it follows for $j=r+1, \ldots, N+1$
\begin{eqnarray*}
\phi(v^j)-v^j &=& \Gammu(v^j) \left( 1 + \frac{\kappa_{\Gamma}-2\|v^j\|}{\|v^j\| + \kappa_{0}}\right) - v^j \ll \Gammu^{\max} \left( 1 + \frac{\kappa_{\Gamma}-2\|v^j\|}{\|v^j\| + \kappa_{0}}\right) - v^j \\
&\ll& \Gammu^{\max}\frac{1}{2k-1} - v^j \stackrel{(\ref{eq:Wahlvonk})}{<} c-v^j. 
\end{eqnarray*}
Altogether with $(\ref{eq:proof:c-vj})$ it follows for all $j=1, \ldots, N+1$
$$ l(y^j)_{l} =\left( (1-t_{j})c+ t_{j}\phi(v^j)\right)_{l} - v^j_{l} \leq c_{l}-v^j_{l}<0. $$
Let $L_l$ denote the $l^{th}$ row of $L$ and  let $W^{m}$ denote the $m^{th}$ column of $W$. From the above consideration it follows $(L_l)^\top \ll 0$ and from $L_1W^1=1 $ we get $W^1>0$. But then it follows $L_1W^1<0$, a contradiction to $L_lW^1=0$ according to equation \eqref{eq:vollstIII}. So $\tau$ is not complete.
\end{proof}

Now we can deduce the following main theorem. 

\begin{theorem}\label{theo:path-konverge} Let $\phi$ be defined as in $(\ref{eq:phix})$ and assume that $\Gamma$ is irreducible and that the operator $\Gammu:\R^N_{+} \rightarrow \R^N_{+}$, deduced from the gain matrix $\Gamma$, satisfies the small gain condition $(\ref{eq:Merrill:SGB})$. Let $\delta>0$ be chosen as in $(\ref{eq:Wahlvondelta})$ with $k \in \N$ according to $(\ref{eq:Wahlvonk})$. Then the simple path starting with $\eta^0$ is finite and the SFP-algorithm converges to a decay point $s \in \Omega(\Gammu)$.\end{theorem}

\begin{proof} The dark grey painted region in Figure
    \ref{fig:posOrthTriang} is compact. The path of complete simplices
    starts in the interior of this region. Theorem \ref{satz:intOrth} and
    Theorem \ref{satz:nichtnachIV} now show that under the above
    assumptions the path starting with $\eta^0$ cannot leave this
    region. So the path remains in this region. Since the region is
    compact there exist only finitely many simplices in this region and we
    are in the situation of Theorem \ref{theo:PathPossibility}~(iv)
    . So the path is finite and ends up in a simplex $\tau \in \R^N_+
    \times \{1\}$ which contains an approximate fixed point of $\phi$ by
    Proposition \ref{satz:Merrill-ApproxFP}. So refining the triangulation
    leads to the convergence of the SFP-algorithm to a fixed point $s$ of
    $\phi$.  Since $\Gammu$ satisfies the small gain condition
    $(\ref{eq:Merrill:SGB})$ it follows by Proposition \ref{satz:FPAP}
    that the fixed point $s$ of $\phi$ lies in the set of decay
    $\Omega(\Gammu)$. So the SFP-algorithm converges to a decay point $s
    \in \Omega(\Gammu)$. \end{proof}

\subsection{Improvement of the algorithm}\label{subsec:ImproveMerrill}

To summarize implementation details we give some suggestions for the choice of $c$, $\delta$, and the constants $\kappa_{0}, \kappa_{h}$ and $\kappa_{\Gamma}$ for a given function $\Gammu$ of dimension $N$.

\paragraph{Suggestions for the choice of $\kappa_{h}, \kappa_{\Gamma}, \kappa_{0}$} \textcolor{white}{$nichts$}\\
Theorem \ref{satz:FPAP} says that a fixed point $s=\phi(s)$ can only lie in region $\1$. Since $\phi(s)=\Gammu(s)+ (\kappa_{h}-2\|s\|)e$ for $s \in \1$ we may expect the fixed point $s$ to have a norm near $\kappa_{h}/2$. So we choose $\kappa_{h}$ as the double size of the norm of the desired fixed point. \\
Several computational experiments have shown that values for
$\kappa_{\Gamma}$ near $\kappa_{h}$ and $\kappa_{0}$ small will probably
lead to small computing times. So we give the suggestions $
\kappa_{\Gamma}=\kappa_{h}+1$ and $\kappa_{0}=1.$ Note that for smaller
values the regions $\2'$ and $\3'$ are small and so the path tends to
leave the region $\1' \cup \2' \cup \3'$ more often. This leads to more
pivoting steps and so to longer computing times.

\paragraph{Suggestion for the choice of $c$} \textcolor{white}{$nichts$}\\
If we have no advance information about the location of the fixed point we
choose $c$ by default as $ c= 0.99\frac{\kappa_{h}}{2\sqrt{N}} [1 \ldots
1]^\top.$
The norm of $c$ then is $\|c\|=0.99\kappa_{h}/2$, thus near $\kappa_{h}/2$ where we expect the fixed point. In addition no direction is preferred.\\
In some cases we have some information about the approximate location of
the decay point. Then we can use this information by using this expected
point as $c$ (if it lies in $\1 \cup \2$) to arrive smaller computing
times.

\paragraph{Suggestion for $\delta$ and the refinement of $\delta$} \textcolor{white}{$nichts$}\\
The choice of $\delta>0$ as in $(\ref{eq:Wahlvondelta})$ is one that leads
to provable convergence but we have seen in experiments that this choice
leads to longer computing times. So we will ignore the choice of $\delta$
as in $(\ref{eq:Wahlvondelta})$ and give another suggestion. To ensure
that the algorithm converges stop the iteration, if the path leaves
the region $\1' \cup \2' \cup \3'$, and start again with the same starting
point and a refined, i.e., smaller $\delta$.
Since the dimension $N$ of the operator $\Gammu$ can get large it is advisable not to choose $\delta$ too small. Since the simplex has a diameter of $\sqrt{N}\delta$ we also have to include the dimension $N$ into the choice of $\delta$ such that the path doesn't leave the region $\1' \cup \2' \cup \3'$. Our suggestion therefore is $ \delta= \frac{\kappa_{h}}{N}.$ \\
The algorithm refines $\delta$ until the desired accuracy is reached. So
it is important that we do not only choose $\delta$ suitably, but even
determine the refining sequence $\{\delta_{k}\}_{k \in \N}$ such that the
approximation is quite good.
Saigal \cite[Section 5]{Sai1977} gave such a sequence and showed that the algorithm converges quadratically, if we assume in addition that $\phi$ is continuously differentiable and its derivative is Lipschitz continuous. \\
In our case we are content with the refining factor $\frac12$, i.e.,
$\delta_{k+1}= \frac{\delta_{k}}{2}$. The simple reason is that the
algorithm stops, if it finds a point in the set of decay. Again numerical
experiments suggest that we do not have to refine often to find a decay
point.

\paragraph{Summary} \textcolor{white}{$nichts$}\\
We want to summarize our suggestions. Let \texttt{norm}$>0$ denote the desired norm of the decay point. Then we have the suggested values to be computed as
$$
\boxed{
\kappa_{h}= 2\texttt{norm},  \qquad \kappa_{\Gamma}=\kappa_{h}+1, \qquad \kappa_{0}=1, \qquad c= 0.99 \frac{\kappa_{h}}{2\sqrt{N}}e, 
\qquad \delta = \frac{\kappa_{h}}{N}.}$$

\medskip

\begin{remark} Note that if the SFP-algorithm does not converge for $\delta$ as in \eqref{eq:Wahlvondelta}, then the small gain condition $\Gammu \not \geq \id$ cannot hold on the whole region $\1\cup\2\cup\3\cup\4$. So we have to choose a smaller \texttt{norm}$>0$.
\end{remark}

%%%%%%%%%%%%%%%%%%%%%%%%%%%%%%%%%%%%%%%%%%%%%%%%%%%%%%%%%%%%%%%%%%%
%%%	EXAMPLES
\section{Examples}\label{sec:Exa}

In this section we give two examples. First note that in \cite{RufferWirth:2010:Stability-verification-for-monotone-syst:} an algorithm is developed to compute decay points, which is derived from a homotopy algorithm due to Eaves \cite{Eaves}. In \cite{RufferDowerIto:2010:Applicable-comparison-principles-in-larg:} an example is given and decay points are computed with the algorithm from \cite{RufferWirth:2010:Stability-verification-for-monotone-syst:}, referred to as Eaves algorithm. 
We state the principle ideas of this article, pick up the results given by Eaves algorithm, and compare them to those of our SFP-algorithm.\\
The second example concerns about a biochemical control circuit model which leads to a nonlinear gain matrix $\Gamma$. We give a general example of monod kinetics and state some conclusions about the input-to-state stability  of this control circuit model. At the end we consider a perturbed system and use the methods presented here to check the local input-to-state stability numerically.

\subsection{Quasi-monotone systems}\label{subsec:QMS}

The motivation for this example was the article of R{\"u}ffer et al. \cite{RufferDowerIto:2010:Applicable-comparison-principles-in-larg:}. Therein a nonlinear system is given and decay points are computed with  
Eaves algorithm from \cite{RufferWirth:2010:Stability-verification-for-monotone-syst:}.
For this purpose a nonnegative matrix  $P \in \R^{N \times N}$ with spectral radius $\rho(P)<1$ is constructed for given dimension $N$. By Perron-Frobenius theory it follows that the matrix $A:=-I_{N}+P$ then has spectral abscissa $\alpha(A):=\max\{\texttt{Re} \lambda : \lambda \text{ is an eigenvalue of } A\} = -1 + \rho(P)<0.$
So
 the matrix $A$ is Hurwitz with negative diagonal entries and nonnegative off-diagonal entries. Now we define a smooth coordinate transformation $S: \R^N \rightarrow \R^N$ by
$$S(v)_{i} = \left\{\begin{array}{ll} e^{v_i-1} & \text{ if }v_i>1 \\v_{i} & \text{ if }v_i \in [-1,1] \\-e^{-v_i-1} & \text{ if }v_i<-1\end{array}\right. .$$
It holds $S(0)=0$ and $S(\R^N_{+}) = \R^N_{+}$. The mapping $S: \R^N_{+} \rightarrow \R^N_{+}$ is a monotone operator. Then the systems
\begin{equation}\label{vSv} \dot v = S'(S^{-1}(v))AS^{-1}(v)=:g(v) \end{equation}
and
\begin{equation}\label{zAz} \dot z = Az=:h(z) \end{equation}
are equivalent under a nonlinear change of coordinates. Let $v^*$ be any decay point for the function $g$ in equation $(\ref{vSv})$.
With it $z^* := S^{-1}(v^*)$ is a decay point for the function $h$ in equation $(\ref{zAz})$. We want to pick up the associated run times and numbers of iterations and compare them with those of the SFP-algorithm.\\

The following  results correspond to matrices $P \in \R^{N \times N}_{+}$ with positive entries in $[0,1]$ generated by a numerical approximation of the uniform distribution, and $30\%$ of those are set to zero. Then $\alpha(A)=-0.2$. The numbers are averages over 100 simulations. Here we assumed $\texttt{norm} =10$, i.e., the norm of the desired decay point $v^*$ is $\| v^* \| \approx 10$. In Table \ref{table:Eaves K1} the results of \cite{RufferDowerIto:2010:Applicable-comparison-principles-in-larg:} are listed. In Table \ref{table:Merrill} we give the results of the SFP-algorithm. In addition, we tested the SFP-algorithm even for large $N$.

\begin{table}[ht] 
%\begin{center} 
%\begin{minipage}[t]{0.4\textwidth}
\begin{center}
\begin{tabular}{c | l | r}
N & run time & \# iterations \\
\hline 
5 & 0.11465s & 267.62  \\
10 & 0.64855s & 2059.65  \\
15 & 1.7833s & 5505.78 \\
25 & 7.987s & 19742.84 
\end{tabular}
%\vspace{2.5em}
\caption{Results of Eaves (K1)-algorithm from \cite{RufferDowerIto:2010:Applicable-comparison-principles-in-larg:} for $\texttt{norm}=10$}
\label{table:Eaves K1}
\end{center}
%\end{minipage}
\end{table}
%\hfill
\begin{table}[ht]
\begin{center}
%\begin{minipage}[t]{0.52\textwidth}
\begin{tabular}{c | r | r | c}
N & run time & \# iterations & simulations \\
\hline 
5 & 0.0277s & 20.9 & 100\\
10 & 0.0415s & 34.5 & 100 \\
15 & 0.0618s & 72.3 & 100\\
25 & 0.1710s & 187.8 & 100 \\
\hline 
50 & 1.180s & 688.4 & 100 \\
100 & 13.22s & 2711.9 & 50\\
150 & 78.35s & 6614.3 &  10 \\
200 & 273.6s & 11243.8 &  10 
\end{tabular}
\caption{Results of the SFP-algorithm for $\texttt{norm}=10$}
\label{table:Merrill}
%\end{minipage}
\end{center}
\end{table}

Note that the run times and iterations can only be compared relatively since the simulations are executed on different computers. Nevertheless, our run times are considerably lower and even for relatively large dimensions we are able to compute decay points in a quite acceptable run time.%\\

\begin{table}[ht] \begin{center} \begin{tabular}{c | r | r | c}
N & run times & \# iterations & simulations \\
\hline 
5 & 0.0451s & 61.6  & 10\\
10 & 0.0680s & 62.5  & 10 \\
15 & 0.0879s & 106.6  & 10\\
25 & 0.2977s & 317.0 & 10 \\
\hline 
50 & 2.141s & 1097.1 & 10 \\
100 & 25.07s & 3991.3  & 10\\
150 & 253.2s & 9214.9  & 10 \\
200 & 542.0s & 16252.1 & 10 
\end{tabular}
\caption{Results of the SFP-algorithm for $\texttt{norm}=1000$}
\label{table:Merrillnorm1000}
\end{center}
\end{table}

In Table \ref{table:Merrillnorm1000} we give run times and iteration
numbers for $\texttt{norm}=1000$.  One can see that we have a relatively
small increase of iteration steps and therefore of run times despite a
quite larger norm.  This is a consequence of the fact that we choose the
size of $\delta$, and with it the mesh size of the starting triangulation
$\delta \sqrt{N}$, in dependency of the norm \texttt{norm} ($\delta
=\frac{ 2\texttt{norm}}{N}$).  That is why the algorithm gets close to the
desired decay point in few steps.

\subsection{A biochemical control circuit model}\label{subsec:BCCM}

We consider the following biochemical control circuit model similar to \cite{Smi95}, 
\begin{align} \label{eqn:bccm}
\dot x_{1}(t) &= g(x_{N}(t))-a_{1}x_{1}(t) +u(t) \notag\\
\dot x_{i}(t) &= x_{i-1}(t)-a_{i}x_{i}(t), \qquad i=2, \ldots, N,\\
x(t) &= [x_{1}(t), \ldots, x_{N}(t)]^\top \in \R^N_{+}, \notag
\end{align}
with $a_{i}>0$ constant for all $i=1, \ldots, N$, $u \in \Li([0,\infty); \R)$ and $g: \R_{+} \rightarrow \R_{+}$ a continuously differentiable function with $g(x)>0$ for all $x>0$. In contrast to \cite{Smi95} we added an external input $u$ and do not assume $g$ to be bounded, but we demand the following assumption, which was introduced in \cite{KJ09}.

\begin{assumption} \label{as:A}
There exist $x_{N}^*>0, K>0$ and $\lambda \in (0,1)$ with $ax_{N}^* = g(x_{N}^*)$ and $a= \prod_{j=1}^Na_{j}$ such that
\begin{equation}\label{eq:(A)}  \frac{K+x_{N}^*}{K+x}x \leq a^{-1}g(x) \leq x_{N}^* + \lambda | x-x_{N}^*| \quad \text{for all } x \geq 0. \end{equation}
\end{assumption}

\begin{remark}\label{rem:MK}
The function 
\begin{equation}\label{eq:MKfunc} g(x)= \frac{b x}{c+x}, \qquad x \geq 0, \quad b, c >0\end{equation} 
models the growth rate of cells or micro-organism and is further known as monod kinetics. Let $a>0$ be arbitrary and  $b, c>0$ such that $b>ac$. Then assumption $\ref{as:A}$ is satisfied. This can easily be seen by setting $x_{N}^* := \frac{b-ac}{a}>0$, $K=c>0$ and $\lambda := \frac{c}{c+x_{N}^*}$.
\end{remark}

Following similar calculations as in \cite{KJ09} we get the following result by setting the gains as
\begin{align*} \gamma_{1,j}(s) &\equiv 0 \quad \text{ for $j \neq N$ and } &\gamma_{1,N}(s) &:= \frac{1}{2} \Big( \ln \big(1+ \theta ( \exp(\sqrt{2s})-1)\big)\Big)^2,\\
\gamma_{i,j}(s) &\equiv 0 \quad \text{ for $j \neq i-1$, } &\gamma_{i,i-1}(s)&:= \frac{1}{2}\Big(\ln \big(1+ \zeta(\exp(\sqrt{2s})-1)\big)\Big)^2 ,\end{align*}
with $\Theta \in \left(\max \{\tfrac{K}{K+x_N^*}, \lambda \}, 1\right)$ and $\zeta \in (1, \Theta^{-\tfrac{1}{N-1}})$.

\begin{theorem}\label{theo:BCCM} \cite[Satz 5.5]{Geiselhart} \label{satz:monodLISS} Consider the system $(\ref{eqn:bccm})$ with $g$ defined as in $(\ref{eq:MKfunc})$ with $b,c>0$. If $b>ac$ with $a= \prod_{j=1}^Na_{j}>0$ then the equilibrium solution $x^*:=[x_{1}^*, \ldots, x_{N}^*]^\top \gg 0$ with
$x_{N}^*:= \frac{b-ac}{a}$ and $x_{i}^*:=\left(  \prod_{j=i}^{N-1}a_{j+1}\right) x_{N}^*$ for $i=1, \ldots, N-1$ is ISS on $\R^N_+\backslash\{0\}$. \end{theorem}

\subsubsection{A perturbed biochemical control circuit model}\label{subsubsec:perturbed}

\begin{figure}[ht]
\centering
\begin{overpic}[scale=.30]{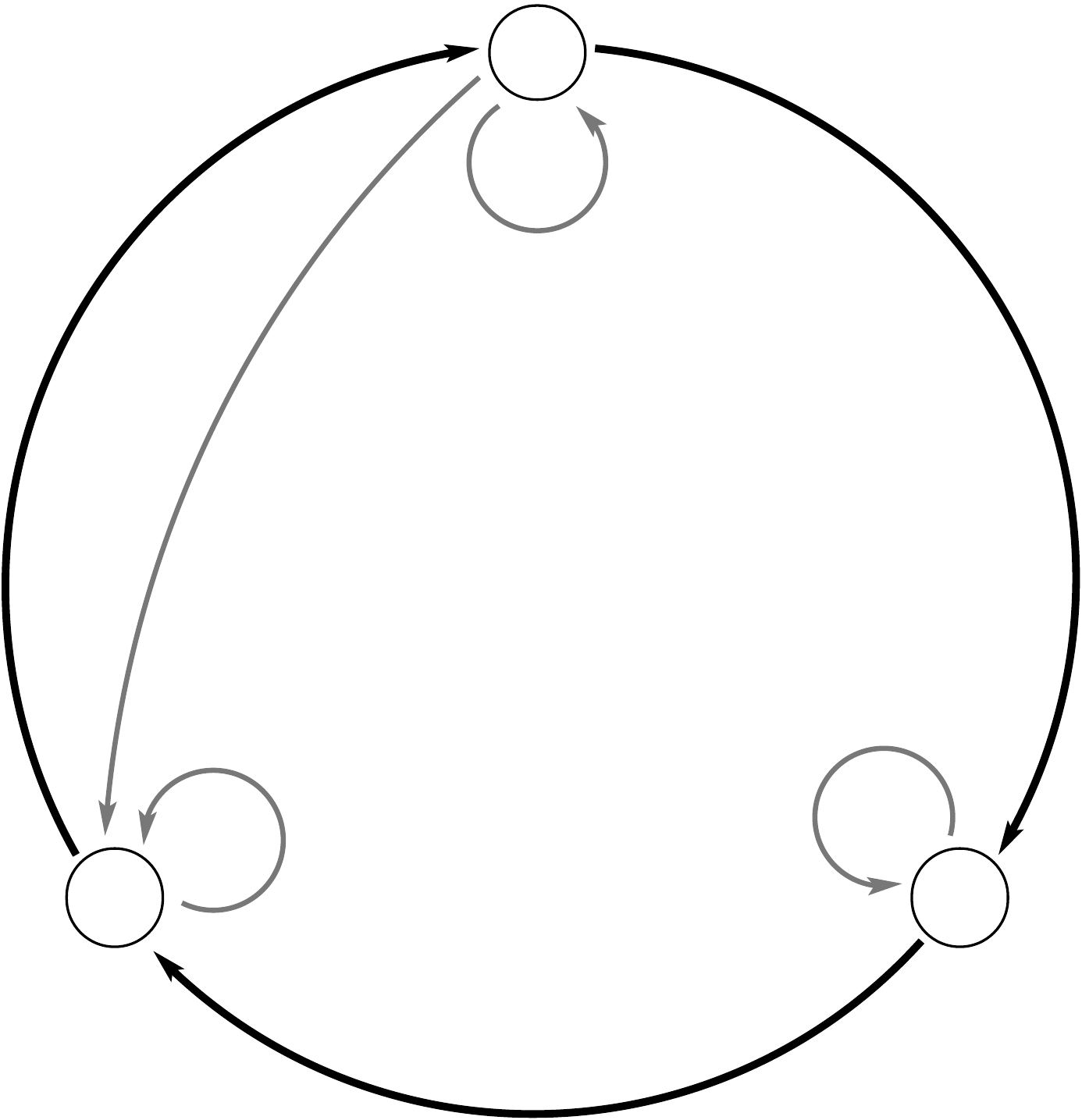} 
\put(46.3, 93.4){${\scriptstyle 1}$} 
\put(84.3,17.7){${\scriptstyle 2}$} 
\put(8.5,17.7){${\scriptstyle 3}$}  
\end{overpic}
\caption{The coupling graph of a perturbed control circuit model} 
\label{fig:GraphBCCM}
\end{figure}

Theorem \ref{satz:monodLISS} is a nice theoretical result. But in applications we are always faced with perturbed systems. In this section we want to check the local input-to-state stability of a perturbed system with the methods developed in this work.\\
For this purpose consider the graph in Figure \ref{fig:GraphBCCM}. The black arcs describe the real couplings of the biochemical control circuit model and the grey arcs describe the perturbations. The underlying system to this coupling graph is given by
\begin{align}
\dot x_{1}(t) & = g(x_{3}(t))-a_{1}x_{1}(t)+u(t) \textcolor{grey}{ + \tilde \varepsilon_{11}(x_{1}(t))} \notag \\
\label{eq:perturbedBCCM} \dot x_{2}(t) & = x_{1}(t)-a_{2}x_{2}(t) \textcolor{grey}{ + \tilde \varepsilon_{22}(x_{2}(t))} \\ 
\dot x_{3}(t) & = x_{2}(t)-a_{3}x_{3}(t) \textcolor{grey}{ + \tilde \varepsilon_{33}(x_{3}(t)) + \tilde \varepsilon_{31}(x_{1}(t))}. \notag
\end{align}

For this system the functions $\tilde \varepsilon$ describe the perturbations. Further let $a_{1}=2$, $a_{2}=1$, $a_{3}=3$ and $a=a_1 \cdot a_2 \cdot a_3 =6$, and $g: \R_{+} \rightarrow \R_{+}$ be the monod function  given in (\ref{eq:MKfunc}) with $b=8$ and $c=1$. 
Then the associated gain matrix $\Gamma$ is of the form
$$ \Gamma = \left(\begin{array}{ccc}
\textcolor{grey}{\tilde \gamma_{11}} & 0 &  \gamma_{13} \\ 
\gamma_{21} & \textcolor{grey}{\tilde \gamma_{22}} & 0 \\ 
\textcolor{grey}{\tilde \gamma_{31}} & \gamma_{32} & \textcolor{grey}{\tilde \gamma_{33}}
\end{array}\right).$$
Here let
\begin{align*}
\gamma_{13}(s) & := \frac{1}{2} \left( \ln(1+ \theta(\exp(\sqrt{2s})-1)) \right)^2 \\
\gamma_{21}(s) & := \frac{1}{2} \left( \ln(1+ \zeta(\exp(\sqrt{2s})-1)) \right)^2 \\
\gamma_{32}(s) & := \frac{1}{2} \left( \ln(1+ \zeta(\exp(\sqrt{2s})-1)) \right)^2 
\end{align*}
be the gain functions from the previous subsection with
\begin{equation}\label{eq:BCCM:thetazeta} \theta \in \left( \max \{ \frac{K}{K+x_N^*}, \lambda \},1 \right) \quad \text{ and } \quad \zeta \in \left(1, \theta^{-1/2} \right). \end{equation}
By Remark \ref{rem:MK} it follows $K=1$, $\lambda=\frac34$ and $x_N^*=1/3$, i.e., $\theta \in (\frac34,1)$. Here we assume
$\theta = 0.8$ and $\zeta =1.1$, where (\ref{eq:BCCM:thetazeta}) is satisfied. Since $b=8>6=ac$ and $a=6>0$, it follows by Theorem \ref{satz:monodLISS} that the equilibrium solution $x^* = \left[\begin{array}{ccc}2 &1 &1/3\end{array}\right]^\top$ 
of the unperturbed system ($\tilde \varepsilon_{11}\equiv \tilde \varepsilon_{22} \equiv \tilde \varepsilon_{31} \equiv \tilde \varepsilon_{33} \equiv 0$ respectively $\tilde \gamma_{11}\equiv \tilde \gamma_{22} \equiv \tilde \gamma_{31} \equiv \tilde \gamma_{33} \equiv 0$) is ISS  on $\R^N_+\backslash\{0\}$.

The perturbed gain functions of the system are given as
$$ \tilde \gamma_{11}(s)  := 0.001s, \quad  \tilde \gamma_{31}(s):= 0.005s^2, \quad \tilde \gamma_{22}(s)  := 0.001s^{0.9}, \quad   \tilde \gamma_{33}(s)  := 0.001s^2 .$$
For the monotone aggregation function $\mu= \sum$ we get the monotone operator $\Gammu$.
Now we apply the SFP-algorithm with the suggestions in section \ref{subsec:ImproveMerrill} and 
$\texttt{norm}=12$ to $\Gammu$ and get the decay point as 
$$w = \left[\begin{array}{c}6.54 \\6.90 \\7.33\end{array}\right] \gg \left[\begin{array}{c}6.527 \\6.886 \\7.325\end{array}\right] \approx \Gammu(w).$$
Since $\lim_{k\rightarrow \infty}\Gammu^k(w)=0$, Theorem \ref{theo:LISS} is applicable, so the perturbed system $(\ref{eq:perturbedBCCM})$ is locally ISS.\\

\begin{figure}[ht]
\centering
\begin{overpic}[scale=0.36]{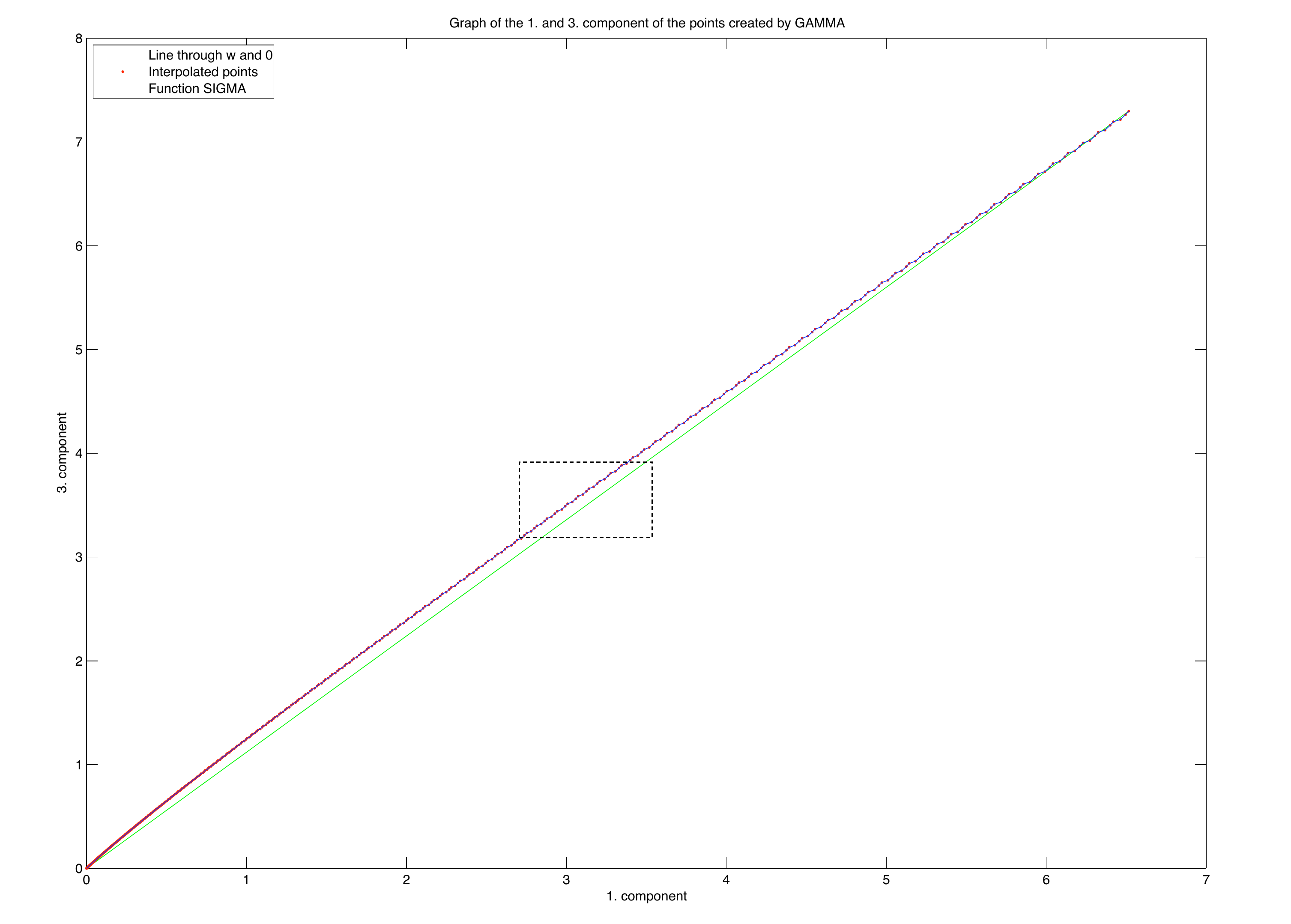}
\put(54,10){\includegraphics[scale=0.128]{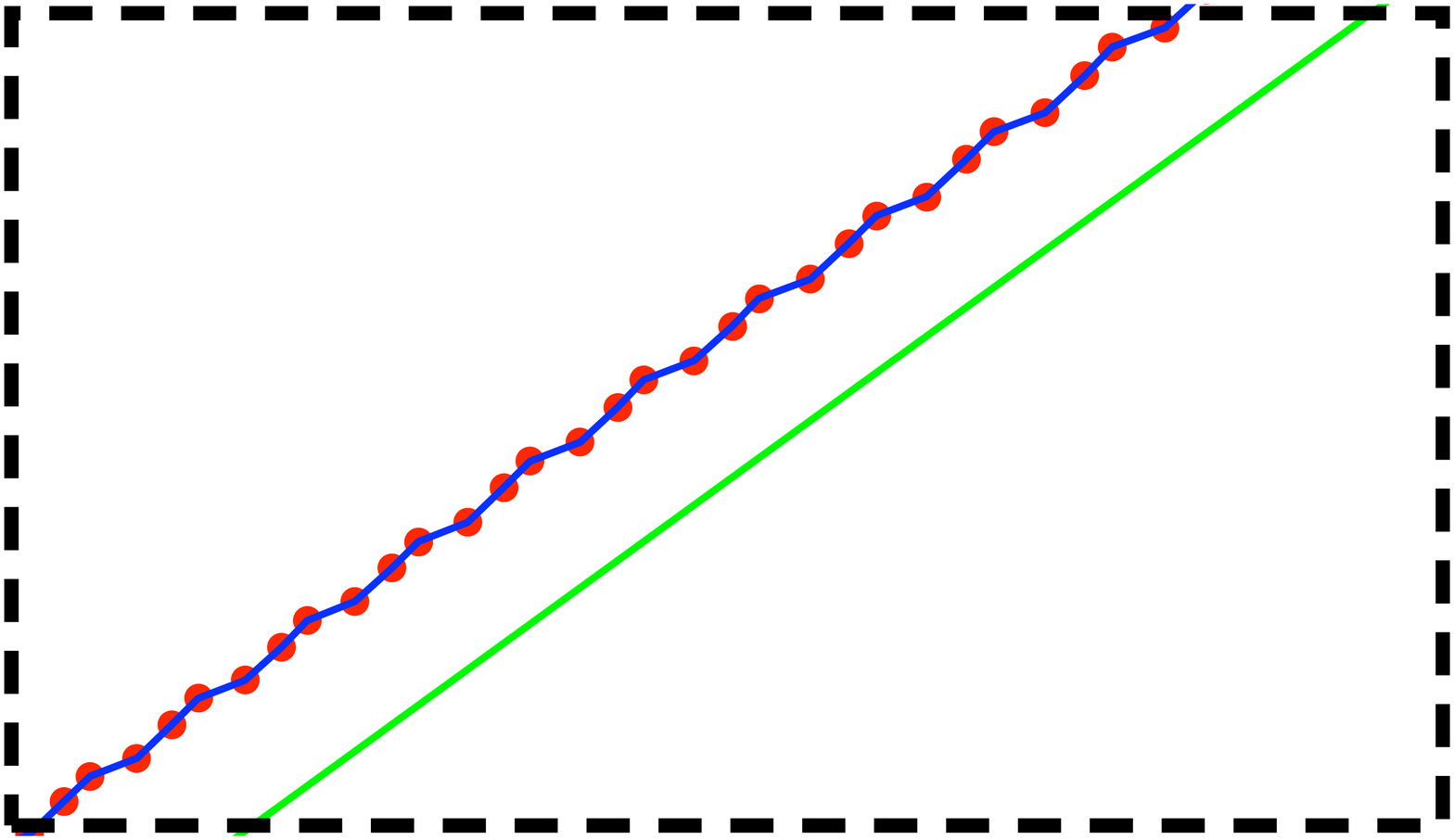}}
\end{overpic}
\caption{Components of the path $\sigma$ from $0$ to $w$} 
\label{fig:sigmaBCCM}
\end{figure}

Finally we illustrate the first and third components of the path $\sigma$ starting in $w$ as defined in \eqref{eq:lininterpol} in Figure \ref{fig:sigmaBCCM}. 
Recall that $\sigma $ is obtained as a linear interpolation (\textcolor{blue}{\pmb{\pmb{---}}}) of the points $\Gammu^k(w)$ (\textcolor{red}{$\bullet$}), where we plot this for $k=1, \ldots, 1000$. 
Indeed, straightforward calculations show that the line from $0$ to $w$ (\textcolor{green}{\pmb{\pmb{---}}}) is not contained in the decay set $\Omega(\Gammu)$. Although the difference between the path and the straight line appears to be negligible we see that without the numerical effort of computing the $\Gammu^k(w)$ no $\Omega$-path is obtained.
% Note that the path $\sigma$ is not linear. For this comparison we additionally plotted the line from $0$ to $w$ (\textcolor{green}{\pmb{\pmb{---}}}). Calculations show in addition that this line does not contain to the set $\Omega(\Gammu)$. 
For a better view we enlarged one region.

\subsubsection{A higher dimensional test}

In this section we consider system (\ref{eqn:bccm}) and want to study the computational effort for higher dimensions. Let $N \in \N$ denote the dimension of system (\ref{eqn:bccm}) with $a_i = \frac{i+1}{i}$, $i=1, \ldots, N$, and the function $g$ be defined as in (\ref{eq:MKfunc}) with $c=1$ and $b =2N$, then Assumption \ref{as:A} is satisfied with
$$a=\prod_{i=1}^N a_i = N+1, \quad x_N^* = \tfrac{N-1}{N+1}, \quad K=1, \quad \lambda = \tfrac{N+1}{2N}.$$
Now by Theorem \ref{theo:BCCM} this system is ISS on $\R^N_+\backslash \{0\}$, if
$$ \theta \in (\lambda,1) = (\tfrac{N+1}{2N},1) \quad \text{and} \quad \zeta \in (1, \theta^{-\frac{1}{N-1}}).$$
In Table \ref{table:BCCM} we tested the computational effort for higher dimensions for \texttt{norm}=12. Note that for large $N$ we have $\theta^{-\frac{1}{N-1}}$ near $1$, so we have only a small range in choosing $\zeta$ such that the system (\ref{eqn:bccm}) still is ISS on $\R^N_+\backslash\{0\}$. So we guess that the decay set will be very thin and so the decay points are harder to reach, resulting in longer computing times. This can be seen in Table \ref{table:BCCM}. Note that the run times include checking that the sequence $\Gammu^k(w)$ is a zero sequence. The counter \texttt{k\_step} indicates the first $\tilde k \in \N$ such that $\|\Gammu^{\tilde k}(w)\|<10^{-9}$.

\begin{table}[ht] \begin{center} \begin{tabular}{c | c | c | r | r | r}
N & $\theta$ & $\zeta$ & run times & \# iterations & \texttt{k\_step} \\
\hline 
10 & 0.75 & 1.020  & 0.30s & 134 & 1215\\
50 & 0.75 & 1.003  & 4.99s & 1405 & 4501\\
70 & 0.75 & 1.002  & 1.72s & 74 & 5911\\
90 & 0.75 & 1.002  & 57.61s & 8426 & 10257\\
\hline
110 & 0.70 & 1.002  & 105.25s & 9632 & 9888\\
150 & 0.70 & 1.001  & 532.43s & 22856 & 8961\\
200 & 0.70 & 1.001  & 2168.18s & 52752 & 12656
\end{tabular}
\caption{Results of the SFP-algorithm for $\texttt{norm}=12$}
\label{table:BCCM}
\end{center}
\end{table}

\section{Conclusions}
\label{sec:conclusions}

In this paper we have presented a homotopy algorithm, that is suitable for
the computation of decay points of gain operators which are crucial in checking the local input-to-state stability. 
The algorithm is proved to converge in a semi-global fashion, 
provided the mesh size of the underlying triangulation is
sufficiently small, but experiments suggest that the result is
conservative and that larger mesh sizes are frequently sufficient. The
algorithm improves on a previous simplicial algorithm. The advantage of
such algorithms is that they can be used to analyze networks with quite
general small gain formulations whereas other approaches rely on special
structure like linearity of the gain operator or the use of maximization
as the monotone aggregation function. In future research we intend to
further develop numerical techniques for small gain results and explore
relevant examples.

%%%%%%%%%%%%%%%%%%%%%%%%%%%%%%%%%%%%%%%%%%%%%%%%%%%%%%%%%%%%%%%%%%%%%%	REFERENCES

\singlespacing %Zeilenabstand verkleinern

\begingroup
%\small
\bibliographystyle{abbrv}
\bibliography{bibliography_NumericalConstruction_7}
\endgroup

\end{document}